\documentclass[12pt]{amsart}
\usepackage{graphicx,url}
\usepackage{epsfig}
\usepackage{fancyhdr}
\usepackage{framed}
\usepackage{floatrow}
\usepackage{mdframed}
\usepackage{amsmath,amsfonts,amssymb,amsthm,mathabx}
\theoremstyle{plain}
\newtheorem{thm}{Theorem}
\newtheorem*{thm*}{Theorem}

\newtheorem{lem}{Lemma}
\newtheorem*{lem*}{Lemma}
\newtheorem{cor}{Corollary}
\newtheorem*{cor*}{Corollary}

\newtheorem{prop}{Proposition}
\newtheorem*{prop*}{Proposition}

\newtheorem{op}{Open Problem}

\theoremstyle{definition}
\newtheorem{example}{Example}
\newtheorem{Def}{Definition}

\newtheorem{Rem}{Remark}

\newcommand{\reftit}{\textit}    
\newcommand{\refis}{\textbf}     
\def\cR{{\mathcal{R}}}
\def\cN{{\mathcal{N}}}

\def\cT{{\mathcal{T}}}

\def\cH{{\mathcal{H}}}
\def\cE{{\mathcal{E}}}
\def\cL{{\mathcal{L}}}
\def\cV{{\mathcal{V}}}

\def\be{\begin{equation}}
\def\ee{\end{equation}}
\def\ben{\begin{equation*}}
\def\een{\end{equation*}}

\begin{document}

\title{Random self-similar trees and\\ a hierarchical branching process}
\author{Yevgeniy Kovchegov}
\address{Department of Mathematics, Oregon State University, Corvallis, OR  97331, USA}
\email{kovchegy@math.oregonstate.edu}

\author{Ilya Zaliapin}
\address{Department of Mathematics and Statistics, University of Nevada, Reno, NV, 89557-0084, USA}
\email{zal@unr.edu}


\subjclass[2000]{Primary 60C05; Secondary 82B99}

\begin{abstract}
We study self-similarity in random binary rooted trees.
In a well-understood case of Galton-Watson trees, a distribution on a space of trees 
is said to be self-similar if it is invariant with respect to the operation of pruning, 
which cuts the tree leaves.
This only happens for the critical Galton-Watson tree (a constant process progeny),
which also exhibits other special symmetries. 
We extend the prune-invariance setup to arbitrary binary trees with edge lengths.
In this general case the class of self-similar processes becomes
much richer and covers a variety of practically important situations.
The main result is construction of the {\it hierarchical branching processes} that 
satisfy various self-similarity definitions (including mean self-similarity and self-similarity in edge-lengths)
depending on the process parameters.
Taking the limit of averaged stochastic dynamics, as the number of trajectories increases, 
we obtain a deterministic system of differential equations that describes the process evolution. 
This system is used to establish a phase transition that separates fading and explosive 
behavior of the average process progeny. 
We describe a class of {\it critical Tokunaga} processes that happen at the phase 
transition boundary.
They enjoy multiple additional symmetries and include the celebrated critical binary 
Galton-Watson tree with independent exponential edge length as a special case. 
Finally, we discuss a duality between trees and continuous functions, and 
introduce a class of {\it extreme-invariant} processes, constructed as the Harris 
paths of a self-similar hierarchical branching process, whose local minima has the 
same (linearly scaled) distribution as the original process. 
\end{abstract}



\date{\today}
\maketitle

\tableofcontents

\section{Introduction}
\label{intro}
Nature commonly exhibits dendritic structures, both static and dynamic, 
that can be represented by tree graphs \cite{Aldous,Vien90,NTG97}.
Examples from diverse applications, together with a review of related 
coalescence and branching models can be found in Aldous \cite{Aldous},
Berestycki \cite{Berestycki}, Bertoin \cite{Bertoin}, Evans \cite{Evans08}, 
Le Gall \cite{LeGall05}, and Pitman \cite{Pitman}.
Despite their apparent diversity, a number of rigorously studied dendritic 
structures possess structural self-similarity, which often allows a low-dimensional 
parameterization \cite{Pec95,NTG97,VG00,KZ15t}.
An illuminating example is the combinatorial structure of river networks, which 
is closely approximated by a two-parametric Tokunaga self-similar model with 
parameters that are independent of river's geographic location
\cite{Tok78,Pec95,DR00,ZZF13}.
Tree self-similarity has been studied primarily in terms of the average values
of selected branch statistics, and rigorous results have been obtained only for 
a very special classes of Markov trees (e.g., binary Galton-Watson trees with no
edge lengths, as in \cite{BWW00}).
At the same time, solid empirical evidence motivates a search for a 
flexible class of self-similar models that
would encompass a variety of observed combinatorial and metric structures 
and rules of tree growth.
We introduce here a general concept of self-similarity
that accounts for both combinatorial and metric tree structure 
(Sec.~\ref{TSSL}, Def.~\ref{def:distss})
and describe a model (Sect.~\ref{HBP}), called {\it hierarchical branching process}, 
that generates a broad range of self-similar trees (Thm.~\ref{main}) and includes the 
critical binary Galton-Watson tree with exponential edge lengths as a special case (Thm.~\ref{main3}).
We study {\it time-invariant} tree distributions, which is
a convenient generalization of Markov growth (Thm.~\ref{main2}).
We also introduce a class of critical self-similar {\it Tokunaga processes} (Sect.~\ref{sec:Tok})
that enjoy additional symmetries --- their edge lengths are i.i.d. random
variables (Prop.~\ref{Tok_tree}), and subtrees of large Tokunaga trees 
reproduce the probabilistic structure of 
the entire random tree space (Props.~\ref{prop:Tok1}).
The duality between planar trees and continuous functions \cite{Harris,Pitman,ZK12} allows us 
using the hierarchical branching process to construct a novel class
of time series that satisfy the {\it extreme-invariance property}: the distribution
of their local minima is the same as that of the original series (Sect.~\ref{LST}).

The paper is organized as follows.
Section~\ref{sec:RSST} introduces the main definitions, including
the Horton-Strahler order of a tree, tree pruning, and a related 
concept of prune-invariance.
Self-similarity for trees with edge lengths is defined in Sect.~\ref{sec:TSS}.
The duality between trees and continuous functions is reviewed in Sect.~\ref{LST}.
In particular, we define here {\it extreme-invariant} processes
that are equivalent to self-similar trees.
The main results are presented in Sect.~\ref{HBP}.
Sect.~\ref{HBP:def} introduces a {\it hierarchical branching process}
that generates a rich collection of self-similar trees.
The hydrodynamic limit for dynamics of the average numbers of Horton-Strahler
branches is established in Sect.~\ref{HBP:hydro}.
The properties of {\it criticality} and {\it time-invariance} are
defined in Sect.~\ref{sec:width} and explored in a self-similar processes
in Sect.~\ref{HBP:ss}.
Critical Galton-Watson process and critical Tokunaga processes, which
generate the most intriguing examples of self-similar trees, are discussed
in Sects.~\ref{sec:cGW}, \ref{sec:Tok}.
Section~\ref{sec:shape} discusses the combinatorial structure of the
critical Tokunaga process.
Section~\ref{open} concludes with two open problems.

\section{Random Trees}
\label{sec:RSST}
The focus of this paper is on finite unlabeled rooted reduced planted binary trees with no planar embedding. 
The space of such trees, which includes the {\it empty tree} $\phi$ comprised of a single 
root vertex and no edges, is denoted by $\cT$.

The existence of the root vertex imposes the parent-offspring relationship  between each pair of the connected vertices in a tree $T\in\cT$: the one closest to the root is called {\it parent}, and the other -- {\it offspring}. The absence of planar embedding in this context means the 
absence of order between the two offspring of the same parent.  
A tree is called {\it reduced} if it has no vertices of degree 2; such trees are also called {\it full} binary trees.  
A tree is called {\it planted} if its root has degree 1.
Accordingly, there are three types of vertices in a tree 
from $\cT \setminus\{\phi\}$: internal vertices of degree 3, 
leaves (degree 1) and the root (degree 1).
The operation of {\it series reduction} removes each degree-two vertex of a 
binary tree by merging its adjacent edges into one.  
Series reduction turns a rooted binary tree into a reduced rooted binary tree.
The edges of a tree from $\cT$ may be assigned positive lengths. 
The space of trees from $\cT$ with edge lengths is denoted by $\cL$.

Any tree from $\cT$ or $\cL$ can be embedded (and represented graphically) in a plane by selecting an order for each pair of offspring of the same parent.
The space of embedded trees from $\cT$ (and respectively $\cL$) is denoted $\cT_{\rm plane}$ (and respectively $\cL_{\rm plane}$).
Examples of trees from $\cL_{\rm plane}$ are found in the bottom row of Fig.~\ref{fig:HST}.
Choosing different embeddings for the same tree $T\in\cT$ (or $T\in\cL$) leads, in general, to different trees from $\cT_{\rm plane}$ (or $\cL_{\rm plane}$).
Sometimes we focus on the combinatorial tree ${\textsc{shape}(T)}\in\cT$, which retains the branching structure of $T$ while omitting its edge lengths and embedding.

\subsection{Tree pruning and related concepts}
\label{pruning}

The concept of self-similarity is related to the pruning operation
\cite{Pec95,BWW00,KZ15t}.
{\it Pruning} (aka Horton pruning) of a tree is an onto function $\cR:\cT\to\cT$,
whose value $\cR(T)$ for a tree $T\ne\phi$ is obtained by removing
the leaves and their parental edges from $T$, followed by series reduction.
We also set $\cR(\phi)=\phi$. 

\begin{figure}[h] 
\centering\includegraphics[width=0.9\textwidth]{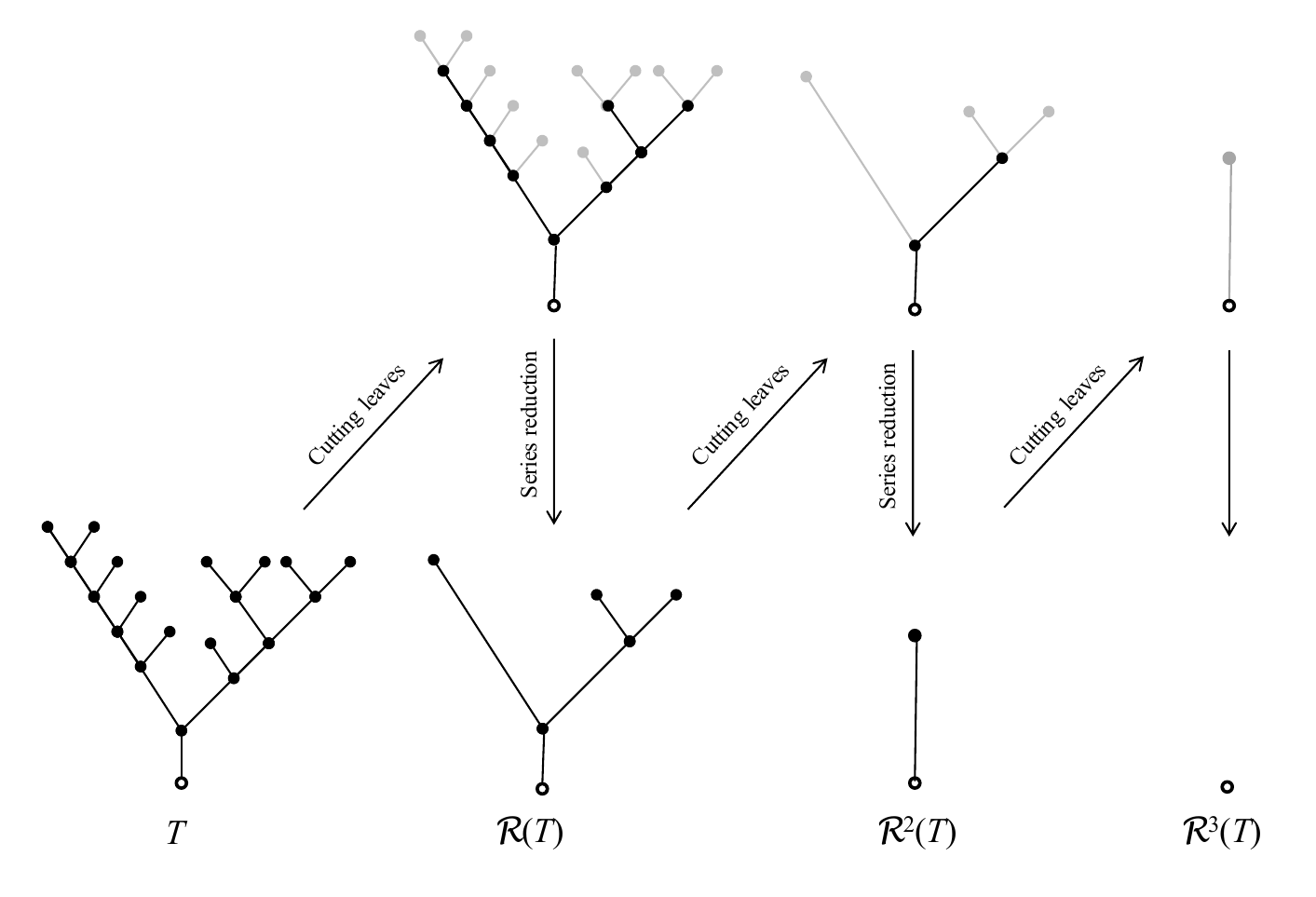}
\caption[Example of Horton-Strahler indexing]
{Example of pruning and Horton-Strahler ordering for a tree $T\in\cL_{\rm plane}$.
The figure shows the two stages of pruning operation --
cutting the leaves (top row), and consecutive series reduction (bottom row).
The initial tree $T$ is shown in the leftmost position of the bottom row.
The gray color in the top row depicts vertices and edges being pruned.
The order of the tree is ${\sf k}(T)=3$, since it is eliminated in three prunings, $\cR^3(T)=\phi$.}
\label{fig:HST}
\end{figure}

\medskip
\noindent
The pruning is also well defined for trees with edge lengths ($\cL$), where series reduction adds the lengths of merging edges,
and for planar trees ($\cT_{\rm plane},\cL_{\rm plane}$), where the embedding of the remaining part of a tree is unaffected by pruning. 
Pruning is illustrated in Fig.~\ref{fig:HST}.

\medskip
\noindent
Pruning induces a contracting map on $\cT$. 
The trajectory of each tree $T$ under $\cR(\cdot)$ is uniquely
determined and finite:
\be\label{TRR0}
T\equiv\cR^0(T)\to \cR^1(T) \to\dots\to\cR^k(T)=\phi,
\ee
with the empty tree $\phi$ as the (only) fixed point.
The pre-image $\cR^{-1}(T)$ of any non-empty tree $T$ consists of an infinite 
collection of trees.
It is natural to think of the distance to $\phi$ under the pruning map
and introduce the respective notion of tree order \cite{Horton45,Strahler}
(see Fig.~\ref{fig:HST}).

\begin{Def}[{{\bf Horton-Strahler order of a tree}}]
The Horton-Strahler order ${\sf k}(T)\in\mathbb{Z}_+$ of a tree $T\in\cT$ is defined as
the minimal number of prunings necessary to eliminate the tree:
\[{\sf k}(T)=\min_{k\ge 0}\left\{\cR^k(T)=\phi\right\}.\]
\end{Def}

\noindent The definition of order is based on the combinatorial shape of a tree.
Accordingly, the order of a tree $T$ from either of spaces 
$\cT_{\rm plane}, \cL,$ or $\cL_{\rm plane}$ is that of $\textsc{shape}(T)$.
The definition implies, in particular, that the order of the empty tree is ${\sf k}(\phi)=0$,
because $\cR^0(\phi)=\phi$.
Most of our discussion will be focused on trees with orders ${\sf k}>0$, and we often
assume that the empty tree $\phi$ has zero probability.

Pruning partitions the tree space $\cT$
into exhaustive and mutually exclusive set of subspaces $\cH_K$ of 
trees of order $K\ge 0$ such that $\cR(\cH_{K+1})=\cH_K$.
Here $\cH_0=\{\phi\}$, $\cH_1$ consists of a single tree comprised of a root and a leaf
connected by an edge, and all other subspaces $\cH_K$, $K\ge 2$, consist of an infinite number of trees.

\begin{figure}[h] 
\centering\includegraphics[width=0.5\textwidth]{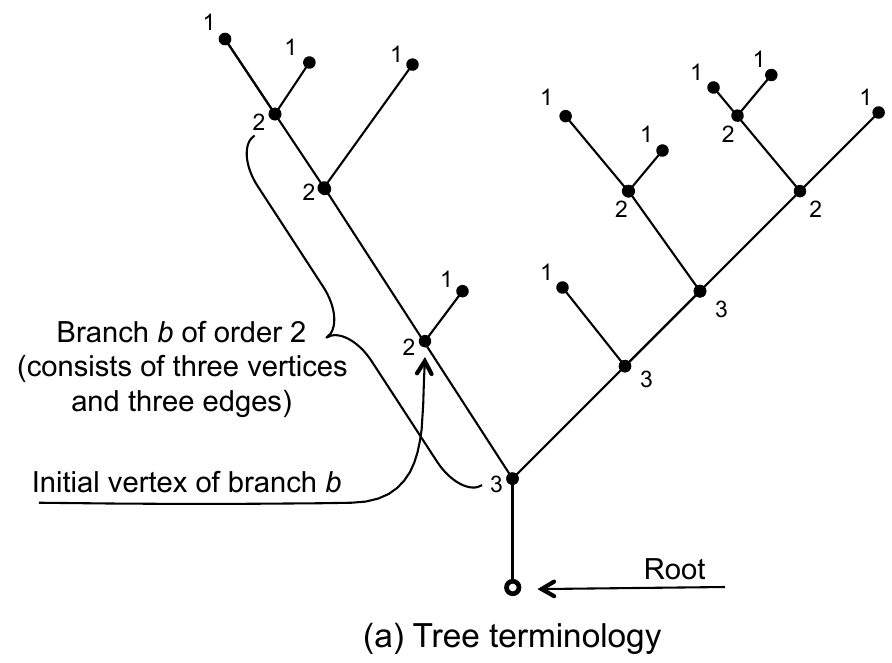}
\centering\includegraphics[width=0.4\textwidth]{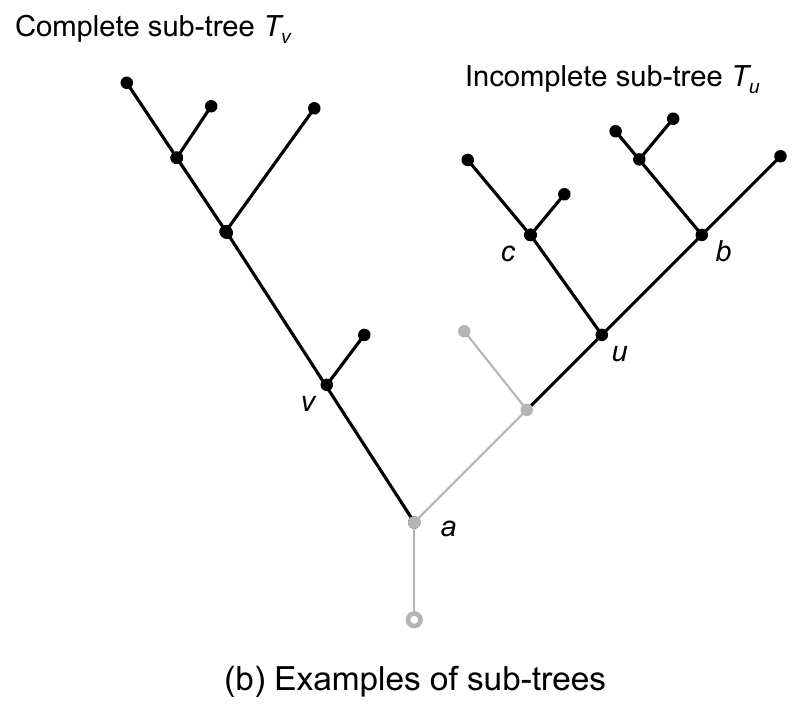}
\caption[Horton-Strahler terminology]
{Illustration of the Horton-Strahler terminology (Def.~\ref{def:HS}) in a 
tree $T$ of order ${\sf k}(T)=3$.
(a) A tree root, a branch, and an initial vertex of a branch.
The numbers indicate the Horton-Strahler orders of the vertices. 
In particular, we have $N_1=10$, $N_2=3$, $N_3=1$, and $N_{1,2}=3$,
$N_{1,3}=1$, $N_{2,3}=1$.
(b) Examples of a complete ($T_v$) and incomplete ($T_u$) subtrees. 
The subtree $T_u$ is incomplete since $u$ is not the initial vertex of a branch. 
The tree shown in this panel has four complete subtrees of order ${\sf k}\ge 2$: $T_v$ of order 2, 
$T_b$ of order 2, $T_c$ of order 2, and $T_a=T$ of order 3.}
\label{fig:terminology}
\end{figure}

\begin{Def}[{{\bf Horton-Strahler terminology}}]
\label{def:HS}
We introduce the following definitions related to the Horton-Strahler order
of a tree (see Fig.~\ref{fig:terminology}):
\begin{enumerate}
\item 
For any non-root vertex $v$ in $T \in \cT \setminus \{\phi\}$, a {\it subtree} $T_v$ of $T$ is defined as the only subtree 
in $T$ rooted at the parental vertex $p(v)$ of $v$, and comprised by $v$ and all 
its descendant vertices together
with their parental edges (Fig.~\ref{fig:terminology}b).
\item
The Horton-Strahler order ${\sf k}(v)$ of a vertex $v\in T$
coincides with the order of the subtree $T_v\in T$ (Fig.~\ref{fig:terminology}a).
\item The parental edge of a vertex has the same order as the vertex.
\item
A connected sequence of vertices of the same order together with their parental edges is called a {\it branch} (Fig.~\ref{fig:terminology}a).
\item
The branch vertex closest to the root is called the {\it initial} vertex of a branch
(Fig.~\ref{fig:terminology}a).
\item
For the initial vertex $v$ of a branch of order $K\le {\sf k}(T)$, a subtree of $T_v$
is called a {\it complete subtree} of order $K$ (Fig.~\ref{fig:terminology}b).
The single complete subtree of order ${\sf k}(T)$ coincides with $T$.
(All subtrees of order ${\sf k}=1$ are complete.)
\end{enumerate} 
\end{Def}

\begin{Rem}
Equivalently, the Horton-Strahler ordering can be done by hierarchical 
counting \cite{DK94,Horton45,Strahler,Pec95,NTG97,BWW00}.
In this approach, each leaf is assigned order ${\sf k}({\rm leaf})=1$.
An internal vertex $p$ whose children have orders $i$ and $j$
is assigned the order 
\[{\sf k}(p)=\max\left(i,j\right)+\delta_{ij} =\lfloor \log_2(2^i+2^j)\rfloor,\] 
where $\delta_{ij}$ is the Kronecker's delta and 
$\lfloor x\rfloor$ denotes the maximal integer less than or equal to $x$.
The Horton-Strahler order of the tree is ${\sf k}(T)=\max\limits_v {\sf k}(v)$, 
where the maximum is taken over all non-root vertices of $T \in \cT \setminus \{\phi\}$.
\end{Rem}

\subsection{Labeling tree vertices}
\label{sec:label}
Sometimes we will need to label the vertices and edges of a tree
(e.g., for selecting a branch or vertex uniformly).
The vertices of a planar tree can be labeled by numbers $1,\dots,\#T$
($\#T$ denoting the total number of vertices in $T$)
in order of depth-first search.
We also assume that label of the parental edge for each vertex is taken
from that vertex.

For a tree with no embedding, labeling is done by selecting a 
suitable embedding and then using the depth-first search labeling as above.
Such embedding should be properly aligned with the pruning
operation, as we describe in the following definition.

\begin{Def}[{{\bf Proper embedding}}]
An embedding function $\textsc{embed}: \cT\to\cT_{\rm plane}$ ($\cL\to\cL_{\rm plane}$) 
is called {\it proper} if for any $T\in\cT$ $(T\in\cL)$
\[\cR\left(\textsc{embed}(T)\right)=\textsc{embed}\left(\cR(T)\right),\]
where the pruning on the left-hand side is in $\cT_{\rm plane}$ ($\cL_{\rm plane}$) and
pruning on the right-hand side is in $\cT$ ($\cL)$.
\end{Def}

A proper embedding for a tree with no edge lengths can be done 
using the following induction construction.
A tree of order ${\sf k}=1$ assumes a unique embedding.
A tree of order ${\sf k}=2$ is embedded by branching all its side-branches
of order 1 to the right.
Assuming there exists a proper embedding for trees of order ${\sf k}\le K$, we construct
the labeling for a tree of order $K+1$. 
All its side-branches (of any order) branch to the right.
To embed the (only) two merging complete subtrees, $\tau_1\ne\tau_2$, of order $K$, 
we consider their farthest non-identical pruning descendants: 
trees $d_i=\cR^{k}(\tau_i)$, $i=1,2$ obtained by the maximal
possible number $k$ of pruning iterations such that $d_1\ne d_2$.
The number $0\le k\le K-2$ is well defined since all trees
of order 1, which is the unltimate pruning limit, coincide.
By construction, the trees $d_i$ differ only by the number of side-branches
of order 1 attached to the tree $d_0=\cR^{k+1}(\tau_i)$, which
already has proper embedding.
Consider the numbers of order-1 side-branches within each edge of $d_0$, in
the order of its labeling:
$(n^{(i)}_1,\dots,n^{(i)}_{\#d_0})$.
The tree whose sequence has the smallest
first non-coinciding number, will branch to the right.

A proper embedding for a tree $T\in\cL$ with edge length
is constructed in the same fashion, with the only correction.
From the two merging complete subtrees of order $K$ with
the same combinatorial structure, the one with the shortest root edge 
branches to the right.
This definition covers the situation of atomless length distribution, which is of primary interest to us.

\section{Tree Self-Similarity}
\label{sec:TSS}
This section defines self-similarity for combinatorial and metric trees.
The term {\it self-similarity} is associated with invariance of a tree distribution 
with respect to the pruning operation $\cR$ introduced in Sect.~\ref{pruning}.
The prune-invariance alone, however, is often insufficient to generate 
interesting families of trees. 
This necessitates an additional property -- {\it coordination} of conditional 
measures on subspaces of trees of a given order.
Coordination together with prune-invariance constitutes the self-similarity
studied in this work. 
 
We start in Sect.~\ref{sec:mss} with a weak form of self-similarity
that only considers the average values of selected branch statistics;
it was introduced in \cite{KZ15t}.
Section~\ref{sec:dss} discusses a stronger version of self-similarity
that operates with tree distributions.
Self-similarity of metric trees is presented in Sect.~\ref{TSSL}.

\subsection{Mean self-similarity of a combinatorial tree}
\label{sec:mss}

Let $\cH_{K}\subset\cT$ be the subspace of trees of Horton-Strahler order $K\ge 0$.
Naturally, $\cH_K\bigcap\cH_{K'}=\emptyset$ if $K\ne K'$, and 
$\bigcup\limits_{K\ge 1} \cH_K =\cT$.
Consider a set of conditional probability measures $\{\mu_{K}\}_{K\ge 0}$ each of which is defined on 
$\cH_K$ by $\mu_K(T) = \mu(T|T\in\cH_K)$ and let $p_K=\mu(\cH_K)$.
Then $\mu$ is represented as a mixture of the conditional measures:
\be
\label{muk}
\mu=\sum_{K=1}^{\infty}p_K\mu_K.
\ee

\noindent We write ${\sf E}_K(\cdot)$ for the mathematical 
expectation with respect to $\mu_{K}$. 
Let $N_k=N_k[T]$ denotes the number of branches of order $k$ in a tree 
$T\in\cT$ (see Fig.~\ref{fig:terminology}a).
We define the {\it average Horton numbers} for subspace $\cH_{K}$ as
\[\cN_{k}[K]= {\sf E}_K(N_k),\quad 1\le k\le K,\quad K\ge 1.\]
Let $N_{i,j}=N_{i,j}[T]$ denote the number of instances when an order-$i$ 
branch merges with an order-$j$ branch, $1\le i <j$, in a tree $T$ (see Fig.~\ref{fig:terminology}a).
Such branches are referred to as {\it side-branches} of order $\{i,j\}$.
Consider the respective expectation $\cN_{i,j}[K]:={\sf E}_K(N_{i,j})$.
The {\it Tokunaga coefficients} $T_{i,j}[K]$ for subspace $\cH_{K}$ are defined as
\be
\label{def_tok}
T_{i,j}[K]=\frac{\cN_{i,j}[K]}{\cN_{j}[K]}, \quad 1\le i <j\le K.
\ee
The Tokunaga coefficient $T_{i,j}[K]$ is hence reflects the average number of side-branches
of order $\{i,j\}$ per branch of order $j$ in a tree of order $K$. 

\medskip
Next, we introduce a property that ensures independence of the side-branch structure 
of a tree order.

\begin{Def}[{{\bf Mean coordination}}]
\label{coord}
A set of measures $\{\mu_K\}_{K\ge 1}$ on $\{\cH_K\}_{K\ge 1}$ is called {\it mean coordinated}
if $T_{i,j}:=T_{i,j}[K]$ for all $K\ge 2$ and $1\le i< j\le K$.
A measure $\mu$ on $\cT$ is called {\it mean coordinated} if the respective conditional measures $\{\mu_K\}$, as in Eq.~\eqref{muk}, are mean coordinated.
\end{Def}
\noindent For a mean coordinated measure $\mu$, 
the Tokunaga matrix $\mathbb{T}_K$ is a $K \times K$ matrix
\[\mathbb{T}_K=\left[\begin{array}{ccccc}
0 & T_{1,2} & T_{1,3} & \hdots & T_{1,K} \\
0 & 0 & T_{2,3} & \hdots & T_{2,K} \\
0 & 0 & \ddots & \ddots & \vdots \\
\vdots & \vdots & \ddots & 0 & T_{K-1,K} \\
0 & 0 & \dots & 0 & 0\end{array}\right],\]
which coincides with the restriction of any larger-order Tokunaga matrix $\mathbb{T}_M$, $M>K$,
to the first $K\times K$ entries.

\begin{Def}[{{\bf Toeplitz property}}]
\label{Tsi}
A set of measures $\{\mu_K\}_{K\ge 1}$ on $\{\cH_K\}_{K\ge 1}$ is said to satisfy 
the {\it Toeplitz property} if $T_{i,j}[K] = T_{j-i}[K]$ for each $K\ge 2$ and some 
sequence $T_k[K]\ge 0$, $k=1,2,\dots$.
The elements of the sequences $T_k[K]$ are also referred to as Tokunaga coefficients, 
which does not create confusion with $T_{i,j}[K]$.
A measure $\mu$ on $\cT$ is said to satisfy the {\it Toeplitz property} if the 
respective conditional measures $\{\mu_K\}$, as in Eq.~\eqref{muk}, satisfy
the Toeplitz property.
\end{Def}

\begin{Def}[{{\bf Mean self-similarity}}]
\label{ss1}
A measure $\mu$ on $\cT$ is called {\it mean self-similar} if it is mean coordinated and satisfies the Toeplitz
property.
\end{Def}

\noindent For a mean self-similar measure the Tokunaga matrix $\mathbb{T}_K$ becomes Toeplitz:
$$\mathbb{T}_K=\left[\begin{array}{ccccc}
0 & T_1 & T_2 & \hdots & T_{K-1} \\
0 & 0 & T_1 & \hdots & T_{K-2} \\
0 & 0 & \ddots & \ddots & \vdots \\
\vdots & \vdots & \ddots & 0 & T_1 \\
0 & 0 & \dots & 0 & 0\end{array}\right].$$

\medskip
Pruning $\cR$ decreases the Horton-Strahler order of each vertex (and hence 
of each branch) by unity; in particular
\be
\label{shift1}
N_k[T] = N_{k-1}\left[\cR(T)\right],\quad k\ge 2,
\ee
\be
\label{shift2}
N_{i,j}[T] = N_{i-1,j-1}\left[\cR(T)\right],\quad 2\le i<j.
\ee
Consider measure $\mu^{\cR}_K$ induced on $\cH_K$ by the pruning operator:
\[\mu^{\cR}_K(A) = \mu_{K+1}\left(\cR^{-1}(A)\right)\quad \forall A\subset \cH_K.\]
The Tokunaga coefficients computed on $\cH_K$ using the induced measure
$\mu^{\cR}_K$ are denoted by $T_{i,j}^{\cR}[K]$.  
Formally,
\be
T_{i,j}^{\cR}[K] = T_{i+1,j+1}[K+1] = \frac{\cN_{i+1,j+1}[K+1]}{\cN_{j+1}[K+1]}.
\ee

\begin{Def}[{{\bf Mean prune-invariance}}]
\label{mpi}
A set of measures $\{\mu_K\}_{K\ge 1}$ on $\{\cH_K\}_{K\ge 1}$ is called 
{\it mean prune-invariant} if $T_{i,j}[K] = T_{i,j}^{\cR}[K]$
(equivalently, $T_{i,j}[K] = T_{i+1,j+1}[K+1]$),
for any $K\ge 2$ and all $1\le i < j\le K$.
A measure $\mu$ on $\cT$ is called {\it mean prune-invariant} if the 
respective conditional measures $\{\mu_K\}$, as in Eq.~\eqref{muk}, are
mean prune-invariant.
\end{Def}

\begin{Def}[{{\bf Mean self-similarity}}]
\label{ss2}
A probability measures $\mu$ on $\cT$ is called 
{\it mean self-similar  with respect to pruning $\cR$} if it is coordinated and mean prune-invariant.
\end{Def}

\begin{prop}
Definitions \ref{ss1} and \ref{ss2} of mean self-similarity are equivalent. 
\end{prop}

\begin{figure}[h] 
\centering\includegraphics[width=\textwidth]{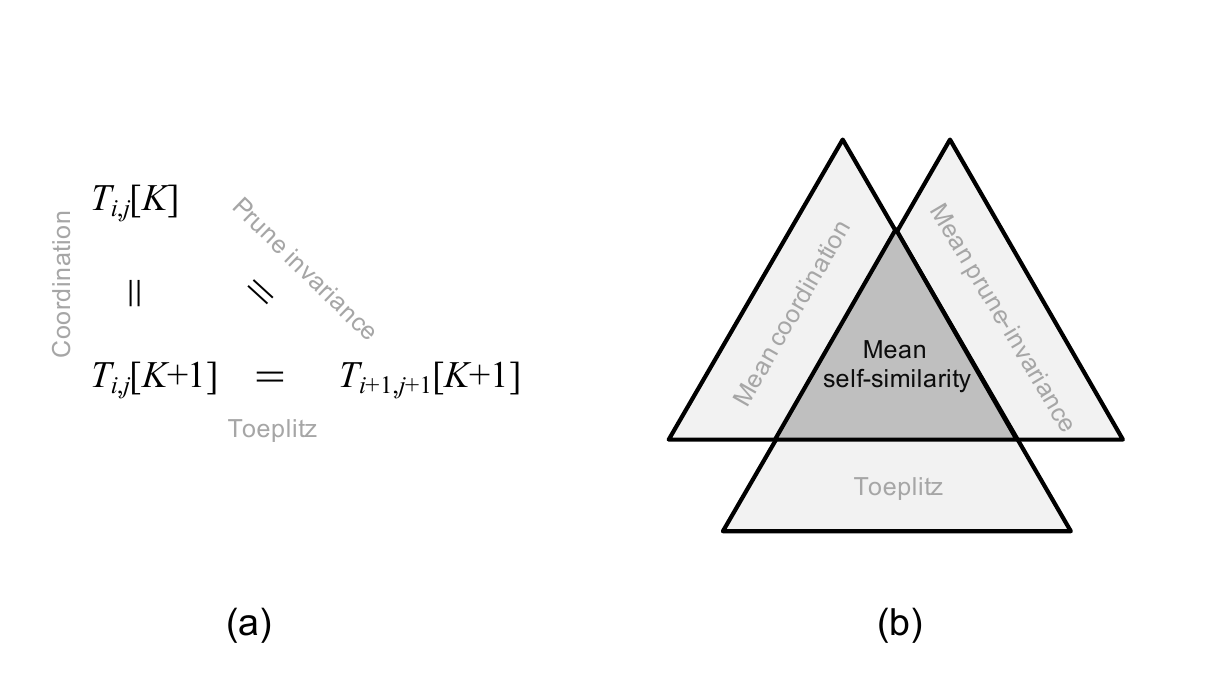}
\caption[Self-similarity]
{Relations among mean coordination, mean prune-invariance, and Toeplitz property.
(a) Pairwise equalities among the quantities $T_{i,j}[K]$, $T_{i,j}[K+1]$, and
$T_{i+1,j+1}[K+1]$ involved in the definitions of mean coordination,
mean prune-invariance, and Toeplitz property.
(b) Venn diagram of the space $\cT$ illustrating the relation
among mean coordination (left triangle), mean prune-invariance (right triangle), 
and Toeplitz property (bottom triangle).
The mean self-similarity (inner dark triangle) is formed by the 
intersection of any pair of the three properties.}
\label{fig:ss}
\end{figure}

\noindent This equivalence was proven in \cite{KZ15t}. 
Its validity is readily seen from the diagram of Fig.~\ref{fig:ss}a, which
shows relations among the quantities $T_{i,j}[K]$, $T_{i,j}[K+1]$, and
$T_{i+1,j+1}[K+1]$ involved in the definitions of coordination,
prune-invariance, and Toeplitz property.
Moreover, we observe that if any two of these properties hold, the third 
also holds. 
The Venn diagram of Fig.~\ref{fig:ss}b illustrates the relation among 
mean coordination, mean prune-invariance, Toeplitz property and mean 
self-similarity in the space $\cT$. 
In this work, we refer to the mean self-similarity with respect to pruning $\cR$
simply as mean self-similarity.

\medskip
A variety of mean self-similar measures can be constructed for 
an arbitrary sequence of Tokunaga coefficients $T_k\ge0$, $k\ge 1$.
Next, we give a natural example \cite{KZ15t}.

\subsection{Example of a mean self-similar measure: Independent random attachment}
The subspace $\cH_1$, which consists of a single-leaf tree,
possesses a trivial unity mass measure.
To construct a random tree from $\cH_2$, we select a discrete probability
distribution $P_{1,2}(n)$, $n=0,1,\dots$, with the mean value $T_1$.
A random tree $T\in\cH_2$ is obtained from the single-leaf 
tree $\tau_1$ of order 1 via the following two operations.
First, we attach two offspring vertices to the leaf of $\tau_1$.
This creates a tree of order 2 with no side-branches -- one 
internal vertex of degree 3, two leaves, and the root.
Second, we draw the number $\tilde N_{1,2}$ from the distribution $P_{1,2}$, and
attach $\tilde N_{1,2}$ vertices to this tree so that they form side-branches 
of order $\{1,2\}$.

In general, to construct a random tree $T\in\cH_{K}$ of order $K\ge 2$ we select a set 
of discrete probability distributions $P_{k,K}(n)$, $k=1,...,K-1$, on $\mathbb{Z}_+$ 
with the respective mean values $T_k$.
A random tree $T\in\cH_{K}$ is constructed by adding branches of order 1 (leaves)
to a random tree $\tau\in\cH_{K-1}$.
First, we add two new child vertices to every leaf of $\tau$ hence producing
a tree $\tilde T$ of order $K$ with no side-branches of order 1.
Second, for each branch $b$ of order $2\le j\le K$ in $\tilde T$ we draw a random number 
$\tilde N_{1,j}(b)$ 
from the distribution $P_{j-1,K}$ and attach $\tilde N_{1,j}(b)$ new child vertices to 
this branch so that they form side-branches of order $\{1,j\}$. 
Each new vertex is attached in a random order with respect to
the existing side-branches.
Specifically, we notice that $m\ge 0$ side-branches attached to a branch of order $j$
are uniquely associated with $m+1$ edges within this branch. 
The attachment of the new $\tilde N_{1,j}(b)$ vertices among the $m+1$ edges
is given by the equiprobable multinomial distribution with $m+1$ categories
and $\tilde N_{1,j}(b)$ trials.

The procedure described above generates a set of measures $\{\mu_K\}_{K\ge 1}$ 
on $\{\cH_K\}_{K\ge 1}$ that are mean coordinated by construction (recall that
the mean values $T_k$ of the distributions $P_{k,K}$ are independent of $K$).

Furthermore, observe that
\[N_{i,j}=\sum_{b_i=1}^{N_j}\tilde N_{1,j-i+1}(b_i),\]
\begin{eqnarray}
\label{condarg}
\cN_{i,j}[K]&=&{\sf E}_K(N_{i,j})={\sf E}_K\left({\sf E}_K(N_{i,j}|N_j)\right)
={\sf E}_K(N_j\,T_{j-i})\\
&=&T_{j-i}\,{\sf E}_K(N_j)=T_{j-i}\,\cN_j[K]\nonumber,
\end{eqnarray}
and hence
$T_{i,j}[K]=\cN_{i,j}[K]/\cN_{j}[K]=T_{j-i}$, so the tree is mean self-similar,
according to Def.~\ref{ss1}.

\begin{Rem}
The properties introduced in this section -- mean coordination,
mean prune-invariance, Toeplitz, and mean self-similarity -- are completely
specified by a set of conditional measures $\{\mu_K\}$, and are independent 
of the randomization probabilities $p_K = \mu(\cH_K)$, see Eq.~\eqref{muk}.
\end{Rem}

\begin{Rem}
\label{rem_coord}
The idea of relating tree mean self-similarity (Def.~\ref{ss2}) to mean prune-invariance (Def.~\ref{mpi}) is quite intuitive (see also \cite{BWW00}).
Much less so is the requirement of mean coordination of conditional measures 
(Def.~\ref{coord}), included in the definition of mean self-similarity.
This requirement is motivated by our goal to bridge the measure-theoretic 
definition of self-similarity via the pruning operation (Def.~\ref{ss2})
to a statistical definition via the branch counting (Def.~\ref{ss1}).
In applications, when a handful of trees of different orders is observed, the coordination assumption
allows one to estimate the Tokunaga coefficients $T_{i,j}$ and make inference regarding
the Toeplitz property; see \cite{Pec95,NTG97,DR00,ZZF13}. 
The absence of coordination, at the same time, 
opens a possibility of having a variety of prune-invariant
measures with no Toeplitz constraint, which are hardly treatable in applications. 
To construct a simplest such measure, let select any tree $T_2$ from the pre-image
of the only tree $T_1\in\cH_1$ of order $K=1$ under the pruning operation: $T_2\in\cR^{-1}(T_1)\in\cH_2$.
In a similar fashion, select any tree $T_{K+1}$ from the pre-image of $T_{K}$ for $K\ge 2$.
This gives us a collection of trees $T_K\in\cH_K$, $K\ge 1$ such that $\cR(T_{K+1}) = T_K$.
Assign the full measure on $\cH_K$ to $T_K$: $\mu_K(T_K)=1$.
By construction, the measures $\{\mu_K\}$ are mean prune-invariant.
They, however, may satisfy neither the mean coordination nor the Toeplitz property.
This construction illustrates how one can produce rather obscure collections of trees 
that are mean prune-invariant, providing a motivation for the
coordination requirement adopted in this work.
\end{Rem}

\subsection{Self-similarity of a combinatorial tree}
\label{sec:dss}
This section introduces a distribution-based approach to self-similarity.

\begin{Def}[{{\bf Prune-invariance}}]\label{def:prune}
Consider a probability measure $\mu$ on $\cT$ such that $\mu(\phi) = 0$.
Let $\nu(T)=\mu \circ \cR^{-1}(T) = \mu \big(\cR^{-1}(T)\big)$. (Note that $\nu(\phi)>0$.)
Measure $\mu$ is called invariant with
respect to the pruning operation (prune-invariant) if for any tree $T\in\cT$ we have 
\be
\label{def:pi}
\nu\left(T~|T\ne\phi\right)=\mu(T).
\ee
\end{Def}

\begin{prop}\label{geom:p}
Let $\mu$ be a prune-invariant measure on $\cT$.
Then the distribution of orders, $p_K=\mu(\cH_K)$, is geometric:
\be
\label{prop:1}
p_K=p\left(1-p\right)^{K-1},\quad K\ge 1,
\ee
where $p=p_1=\mu(\cH_1)$,
and for any $T\in\cH_K$
\be
\label{prop:2}
\mu_{K+1}\left(\cR^{-1}(T)\right) = \mu_K(T).
\ee
\end{prop}
\begin{proof}
Pruning $\cR$ is a shift operator on the sequence of subspaces $\{\cH_k\}$:
\be
\label{H_shift}
\cR^{-1}(\cH_{K-1})=\cH_K,~K\ge 2.
\ee
The only tree eliminated by pruning is the tree of order 1: 
$\{\tau:\cR(\tau)=\phi\}=\cH_1.$
This allows to rewrite \eqref{def:pi} for any $T\ne\phi$ as
\be
\label{def:pi1}
\mu\left(\cR^{-1}(T)\right)=\mu(T)\left(1-\mu(\cH_1)\right).
\ee
Combining \eqref{H_shift} and \eqref{def:pi1} we find for any $K\ge 2$
\be
\label{mu_shift}
\mu\left(\cH_K\right)
\stackrel{{\rm by}~\eqref{H_shift}}{=}
\mu\left(\cR^{-1}(\cH_{K-1})\right)
\stackrel{{\rm by}~\eqref{def:pi1}}{=}
\left(1-\mu(\cH_1)\right)\mu(\cH_{K-1}),
\ee
which establishes \eqref{prop:1}.
Next, for any tree $T\in\cH_K$ we have
\[\mu(T) = \mu(\cH_1)\left(1-\mu(\cH_1)\right)^{K-1}\mu_K(T),\]
\[\mu\left(\cR^{-1}(T)\right) = \mu(\cH_1)\left(1-\mu(\cH_1)\right)^K\mu_{K+1}\left(\cR^{-1}(T)\right).\]
Together with \eqref{def:pi1} this implies \eqref{prop:2}.
\end{proof}

Proposition~\ref{geom:p} shows that a prune-invariant measure $\mu$ is 
completely specified by its conditional measures $\mu_K$ and the mass 
$p=\mu(\cH_1)$ of the tree of order $K=1$.
The same result was obtained for Galton-Watson trees in \cite[Thm.~3.5]{BWW00}.

\medskip
Next, we introduce a distributional analog of the mean coordination property;
see Def.~\ref{coord} and Remark~\ref{rem_coord}.
Specifically, we assume that a complete subtree $T_K$ of a given order $K$
randomly selected from a random tree $T_H$ of order $H\ge K$ has 
a common distribution independent of $H$.
Since a tree $T_K$ of order $K$ has only one complete subtree of order $K$, which
coincides with $T_K$, this common distribution must be $\mu_K$.
Formally, consider the following process of selecting a {\it uniform random 
complete subtree} ${\sf subtree}_{K,H}$ of order $K$ from a random tree 
$T_H\in\cH_H$.
First, select a random tree $T_H$ according to the conditional measure $\mu_H$. 
Label all complete subtrees of order $K$ in $T_H$ in order
of proper labeling of Sect.~\ref{sec:label}, and
select a uniform random subtree, which we denote ${\sf subtree}_{K,H}$.
By construction,  ${\sf subtree}_{K,H}\in\cH_K$;
we denote the corresponding sampling measure on $\cH_K$ by $\mu^H_K$.

\begin{Def}[{{\bf Coordination}}]
\label{def:coord}
A set of measures $\{\mu_K\}_{K\ge 1}$ on $\{\cH_K\}_{K\ge 1}$ is called
{\it coordinated} if $\mu^H_K(T) = \mu_K(T)$ for any $K\ge 1$, $H\ge K$, and $T\in\cH_K$.
A measure $\mu$ on $\cT$ is called {\it coordinated} if the 
respective conditional measures $\{\mu_K\}$, as in Eq.~\eqref{muk}, are
coordinated.
\end{Def}

\begin{example}
\label{ex1}
The space of finite binary Galton-Watson trees has the 
coordination property. 
Recall that a random binary Galton-Watson tree starts with a single progenitor (root) and increases its depth in discrete steps:
at every step each existing vertex can either split in two with probability $p_2$ or become a leaf (disappear) with 
probability $p_0=1-p_2$.
We denote the resulting tree distribution on $\cT$ by $\mathcal{GW}(p_0,p_2)$. 
This Markovian generation mechanism creates complete subtrees of the same structure, 
independently of the tree order. 
This implies coordination. 
\end{example}

\begin{Def}[{{\bf Combinatorial self-similarity}}]
\label{def:ss}
A probability measure $\mu$ on $\cT$ is called {\it (combinatorially) self-similar with respect to pruning $\cR$} if it is coordinated and prune-invariant.
\end{Def}

In this work, we refer to combinatorial self-similarity with respect to pruning $\cR$
simply as combinatorial self-similarity.
It was established in \cite{BWW00} that critical binary Galton-Watson trees 
$\mathcal{GW}(1/2,1/2)$ are prune-invariant.
Together with coordination (see Example~\ref{ex1}), this implies combinatorial self-similarity.
It also has been shown in \cite{BWW00} that non-critical binary Galton-Watson trees
($p_0\ne 1/2$) are not prune-invariant. 
This gives an example of coordinated measures that are not prune-invariant.
Prune-invariant measures with no coordination can be easily constructed following the 
approach of Remark~\ref{rem_coord}.
To make that construction consistent with the definition of distributional prune-invariance (Def.~\ref{def:prune}), each tree $T_K\in\cH_K$ must be assigned the probability 
$p_K=p(1-p)^{K-1}$.

\subsection{Horton law in self-similar trees}
We say that a random tree $T$ satisfies a {\it strong Horton law} 
if the respective sequence $\cN_k[K]$ of branch numbers decays in geometric fashion as
$k$ increases.
Formally, we require
\be
\label{strongHorton}
\lim_{K\to\infty}\frac{\cN_k[K]}{\cN_1[K]}=R^{1-k},\text{ for any }k\ge 1.
\ee
Horton law and its ramifications, which epitomize scale-invariance of dendritic hierarchical structures, are indispensable in hydrology (e.g., \cite{Shreve66,Pec95,DR00}) and have been reported in biology and other areas; see \cite{NTG97,KZ18} and references therein. 
It has been shown in \cite{KZ16} that the tree that describes a trajectory Kingman's coalescent process with $N$ particles obeys a weaker version of Horton law as $N\to\infty$, and that the first pruning of this tree for any finite $N$
is equivalent to a level set tree of a white noise (see Sect.~\ref{LST} for definitions).

A necessary and sufficient condition for the strong Horton law in a mean self-similar tree
has been established in \cite{KZ15t}:
$$\limsup_{k\to\infty} T_k^{1/k}<\infty.$$
The Horton exponent $R$ in this case is given by $R=1/w_0$, where $w_0$ is the
only real root of 
\[\hat t(z) = -1 + 2z +\sum_{k=1}^{\infty} z^k\,T_k\]
within the interval $(0,1/2]$.
Informally, this means that any mean self-similar tree with a ``tamed''
sequence of Tokunaga coefficients satisfies the strong Horton law.

\subsection{Self-similarity of a tree with edge lengths}
\label{TSSL}
Consider a tree $T\in\cL$ with edge lengths given
by a positive vector $l_T=(l_1,\dots,l_{\#T})$ and let
$\textsc{length}(T)=\sum_il_i$.
We assume that the edges 
are labeled in a proper way as described in Sect.~\ref{sec:label}.
A tree is completely specified by its combinatorial shape
$\textsc{shape}(T)$ and edge length vector $l_T$.
The edge length vector $l_T$ can be specified by
distribution $\chi(\cdot)$ of a point $x_T=(x_1,\dots,x_{\#T})$ on the simplex
$\sum_i x_i = 1$, $0<x_i\le 1$, and conditional distribution $F(\cdot|x_T)$ of the tree length
$\textsc{length}(T)$, where
\[l_T = x_T\times\textsc{length}(T).\]
A measure $\eta$ on $\cL$ is a joint distribution of 
tree's combinatorial shape and its edge lengths; it has
the following component measures. 
\begin{align*}
&{\rm Combinatorial~shape:}\quad \mu(\tau)={\sf Law }\left(\textsc{shape}(T)=\tau\right),\\
&{\rm Relative~edge~lengths:}\quad \chi_\tau(\bar x) ={\sf Law }\left(x_T=\bar{x} \,|\,\textsc{shape}(T)=\tau \right),\\
&{\rm Total~tree~length:}\quad F_{\tau, \bar x} (\ell) ={\sf Law }\left(\textsc{length}(T)=\ell \,|\,x_T=\bar{x}, ~ \textsc{shape}(T)=\tau \right). 
\end{align*}

\noindent The definition of self-similarity for a tree with edge lengths builds
on its analog for combinatorial trees in Sect.~\ref{sec:dss}.
The combinatorial notions of coordination (Def.~\ref{def:coord}) 
and prune-invariance (Def.~\ref{def:prune}), which we refer to as coordination 
and prune-invariance {\it in shapes}, are complemented with analogous properties 
{\it in edge lengths}. 
Formally, we denote by $\mu^H_K(\tau)$, $\chi^H_\tau(\bar x)$, and $F^H_{\tau, \bar x}(\ell)$ 
the component measures for a uniform complete subtree ${\sf subtree}_{K,H}$.
(Notice that the subtree order $K$ is completely specified by the tree shape $\tau$,
which explains the absence of subscript $K$ in the component measures for subtree length).
We also consider the distribution of edge lengths after pruning:
$$\Xi_\tau(\bar x) ={\sf Law }\left(x_{\cR(T)}=\bar{x}  \,|\,\textsc{shape}\big(\cR(T)\big)=\tau \right)$$
and
$$\Phi_{\tau, \bar x} (\ell)={\sf Law }\left(\textsc{length}\big(\cR(T)\big)=\ell \,|\,x_{\cR(T)}=\bar{x}, ~\textsc{shape}\big(\cR(T)\big)=\tau \right).$$
Finally, we adopt here the notation $\cH_K$ for a subspace of trees of order $K\ge 1$ from $\cL$, and 
consider conditional measures $\mu_K(\tau)=\mu(\tau|{\sf k}(\tau)=K)$, $K\ge 1$, 
for a tree $\tau\in\cL$.

\begin{Def}[{{\bf Self-similarity of a tree with edge lengths}}]\label{def:distss}
We call a measure $\eta$ on $\cL$ {\it self-similar with respect to pruning $\cR$} if the following conditions hold 
\begin{itemize}
\item[(i)] The measure is coordinated in shapes. This means that
for every $K\ge 1$ and every $H\ge K$ we have
\[\mu^H_K(\tau)=\mu_K(\tau)\qquad \forall \tau \in \cH_K .\] 
\item[(ii)] The measure is coordinated in lengths. This means that
for every $K\ge 1$, $H\ge K$, and $\tau \in \cH_K$ we have
\[\chi^H_\tau(\bar x)=\chi_{\tau}(\bar{x}) \quad \bar{x}\text{-}a.s. ~,\]  
and for every given $\bar x$,
\[ F^H_{\tau, \bar x}(\ell)= F_{\tau, \bar x} (\ell)  \quad \ell\text{-}a.s. \] 
\item[(iii)]
The measure is prune-invariant in shapes.
This means that for $\nu=\mu\circ\cR^{-1}$ we have
\[\mu(\tau)=\nu(\tau|\tau\ne\phi).\]
\item[(iv)]
The measure is prune-invariant in lengths.
This means that 
\[\Xi_\tau(\bar x)=\chi_\tau(\bar x)\]
and there exists a {\it scaling exponent} $\zeta>0$ such that for any 
combinatorial tree $\tau\in\cT$ we have 
\[\Phi_{\tau, \bar x} (\ell)=\zeta^{-1}F_{\tau, \bar x} \left(\frac{\ell}{\zeta}\right).\]
\end{itemize}
\end{Def}
 
In this work, we refer to self-similarity with respect to pruning $\cR$ simply as self-similarity.
Section~\ref{HBP} below introduces a rich class of measures that satisfy 
this definition.

\section{Tree Representation of Continuous Functions}
\label{LST}

We review here the results of \cite{LeGall93,NP,Pitman,ZK12} on tree 
representation of continuous functions.
This allows us to apply the self-similarity concepts to time series 
and motivates discussion in Sects.~\ref{sec:cGW},~\ref{sec:Tok} below.

\subsection{Harris path}
\label{sec:Harris}
\begin{figure}[h] 
\centering\includegraphics[width=0.7\textwidth]{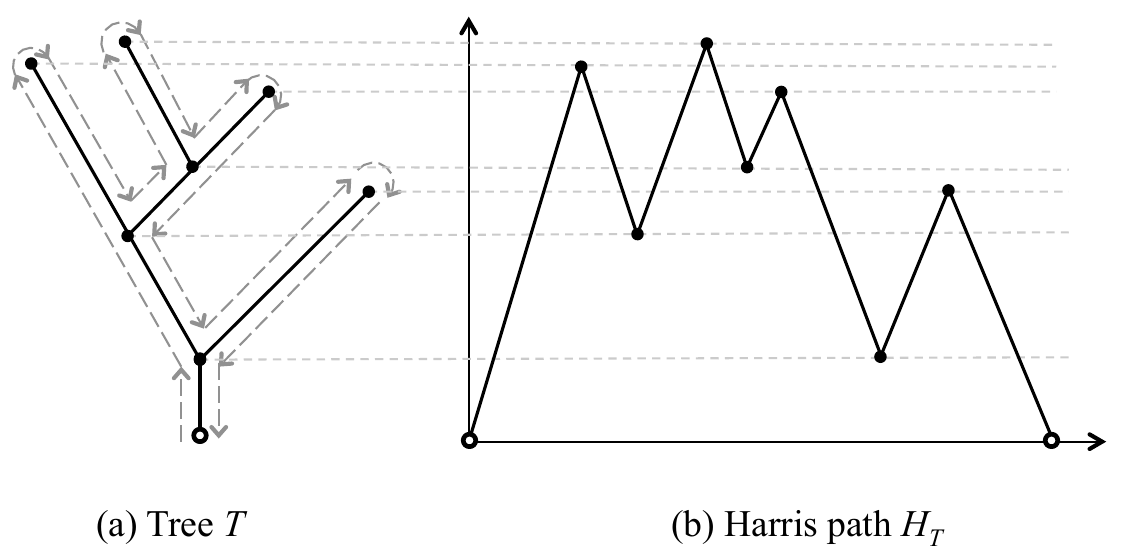}
\caption[Harris path]
{(a) Tree $T$ and its depth-first search illustrated by dashed arrows.
(b) Harris path $H_T(t)$ for the tree $T$ of panel (a). }
\label{fig:Harris}
\end{figure}

For any embedded tree  $T\in\cL_{\rm plane}$ with edge lengths, the {\it Harris path}
is defined as a piece-wise linear function \cite{Harris,Pitman} 
\[H_T(t)\,:\,[0,2\cdot\textsc{length}(T)]\to\mathbb{R}\]
that equals the distance from the root traveled along the tree $T$ 
in the depth-first search, as illustrated in Fig.~\ref{fig:Harris}.
For a tree $T$ with $n$ leaves, the Harris path 
$H_T(t)$ is a piece-wise linear positive excursion
that consists of $2n$ linear segments with alternating slopes $\pm 1$ \cite{Pitman}.

\subsection{Level set tree}
\label{sec:level}
Consider a continuous function $X_t$, $t\in[a,b]$ with a finite number of
distinct local minima.
The level set $\mathcal{L}_{\alpha}\left(X_t\right)$ is defined 
as the pre-image of the function values above $\alpha$: 
\[\mathcal{L}_{\alpha}\left(X_t\right) = \{t\,:\,X_t\ge\alpha\}.\]
The level set $\mathcal{L}_{\alpha}$ for each $\alpha$ is
a union of non-overlapping intervals; we write 
$|\mathcal{L}_{\alpha}|$ for their number.
Notice that 
$|\mathcal{L}_{\alpha}| = |\mathcal{L}_{\beta}|$ 
as soon as the interval $[\alpha,\,\beta]$ does not contain a value of
local extrema of $X_t$ and  
$0\le |\mathcal{L}_{\alpha}| \le n$, where $n$ is the number 
of the local maxima of $X_t$. 

\begin{figure}[h] 
\centering\includegraphics[width=0.6\textwidth]{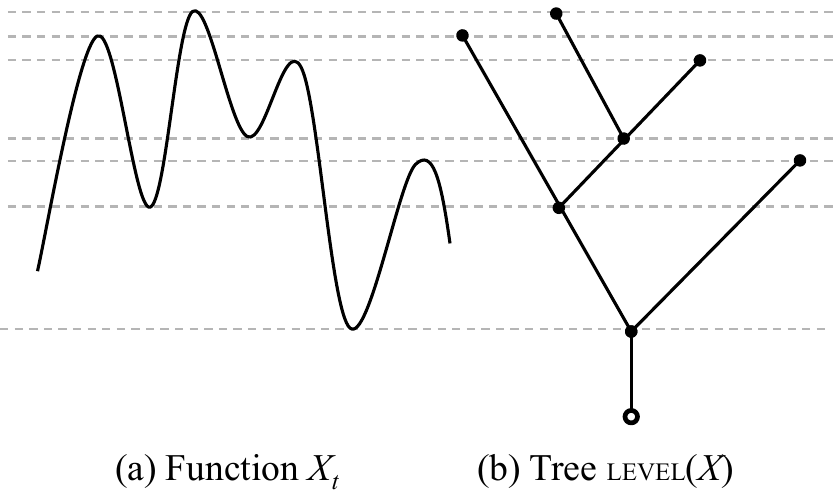}
\caption{Function $X_t$ (panel a) with a finite number of local
extrema and its level-set tree $\textsc{level}(X)$ (panel b).
}
\label{fig:LST}
\end{figure}

The {\it level set tree} $\textsc{level}(X_t)\in\cL_{\rm plane}$ is a 
tree that describes the topology of the level sets $\mathcal{L}_{\alpha}$ 
as a function of threshold $\alpha$, as illustrated in Fig.~\ref{fig:LST}.
Specifically, there are bijections between 
(i) the leaves of $\textsc{level}(X_t)$ and the local
maxima of $X_t$,
(ii) the internal (parental) vertices of $\textsc{level}(X_t)$ 
and the local minima of $X_t$ (excluding possible local minima
at the boundary points), and
(iii) the pair of subtrees of $\textsc{level}(X_t)$ rooted at the 
parental vertex that corresponds to a local
minima $X_{t^*}$ and the first positive excursions (or meanders bounded
by $t=a$ or $t=b$) of $X_t-X_{t^*}$ to right and left of $t^*$.
Every edge in the tree is assigned a length equal the difference
of the values of $X_t$ at the local extrema that correspond to the
vertices adjacent to this edge according to the bijections (i) and (ii) above.
The lowest internal local minimum (achieved for $a<t<b$) corresponds to the first 
descendant of the root, see Fig.~\ref{fig:LST}.
If the global minimum is achieved on the interval boundary ($t=a$ or $t=b$), then it corresponds
to the tree root; otherwise, we artificially add a root to the tree, with an arbitrary root edge length;
this situation is illustrated in Fig.~\ref{fig:LST}.

By construction, the level-set tree $\textsc{level}(X_t)$ is completely
determined by the sequence of the values of local extrema of $X_t$,
and is independent of timing of those extrema, as soon as their order
is preserved.
This means, for instance, that if $g(t)$ is a continuous and monotone increasing 
function on $[a,b]$, then the following trees are equivalent:
\[\textsc{level}(X_t)
=\textsc{level}\left(X_{g(t)}\right).\]
Hence, without loss of generality we can focus on the level set trees
of continuous functions with alternating slopes $\pm 1$.
By $\cE^{\rm ex}$ we denote the space of all positive piece-wise linear continuous finite excursions with alternating slopes $\pm 1$.


By construction, the level set tree of an excursion from $\cE^{\rm ex}$ and 
Harris path are reciprocal to each other as described in the following statement.

\begin{prop}[{{\bf Reciprocity of Harris path and level set tree}}]
The Harris path
$H: \cL_{\rm plane}\to \cE^{\rm ex}$ and the level set tree
${\textsc{level}}: \cE^{\rm ex}\to\cL_{\rm plane}$ are reciprocal
to each other.
This means that for any $T\in\cL_{\rm plane}$ we have
$\textsc{level}(H_T(t))\equiv T,$
and for any $Y_t\in\cE^{\rm ex}$ we have
$H_{\textsc{level}(Y_t)}(t)\equiv Y_t.$
\end{prop}

\begin{figure}[h] 
\centering\includegraphics[width=0.8\textwidth]{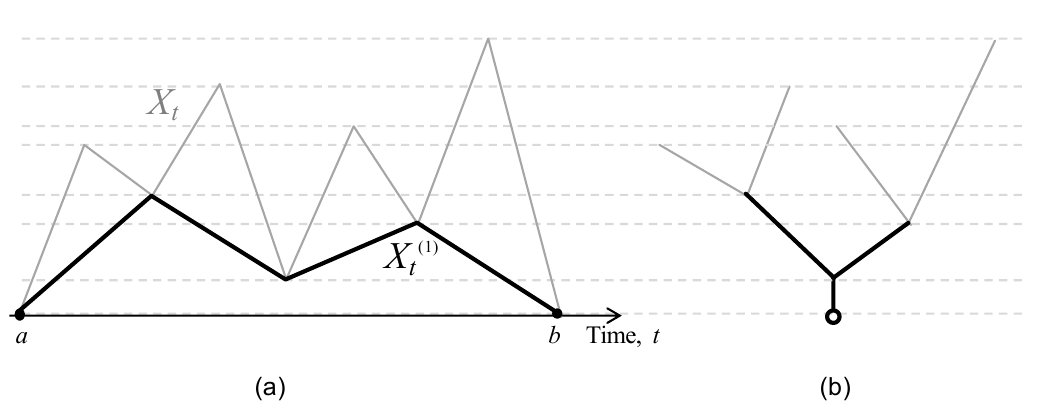}
\caption{Pruning of a positive excursion: transition to the local minima of an excursion 
$X_t$ corresponds to pruning of the corresponding level set tree.
(a) An original excursion $X_t$ (gray line) and linearly interpolated sequence 
$X^{(1)}_t$ of the respective local minima (black line).
(b) The level set tree $\textsc{level}(X^{(1)}_t)$  of the local minima sequence 
(black lines) is obtained by pruning of the level set tree $\textsc{level}(X_t)$
of the original excursion (whole tree). The pruned edges are shown in gray -- each of them
corresponds to a local maximum of the original excursion.}
\label{fig:prune1}
\end{figure}

\subsection{Pruning for positive excursions}
This section examines level set tree and the respective pruning for a positive continuous 
excursion $X_t$ on $[a,b]$. 
Formally, consider a continuous function $X_t$ with a finite number
of local minima and such that $X_a=X_b=0$ and $X_t>0$ for $a<t<b$.
Furthermore, consider excursion $X^{(1)}_t$, $t\in[a,b]$, obtained by a linear interpolation
of the boundary values and the local minima of $X_t$; as well as
functions $X^{(m)}_t$, $t\in[a,b]$, for $m\ge 1$, obtained by taking the local minima of $X_t$ 
iteratively $m$ times, and linearly interpolating their values together with
$X_a=X_b=0$ (see Fig.~\ref{fig:prune1}a).

In the space of level set trees of continuous functions, the pruning $\cR$ corresponds to 
coarsening the respective function by removing (smoothing) its local maxima, 
as illustrated in Fig.~\ref{fig:prune1}. 
An iterative pruning corresponds to iterative transition to the local minima.

\begin{figure}[h] 
\centering\includegraphics[width=0.9\textwidth]{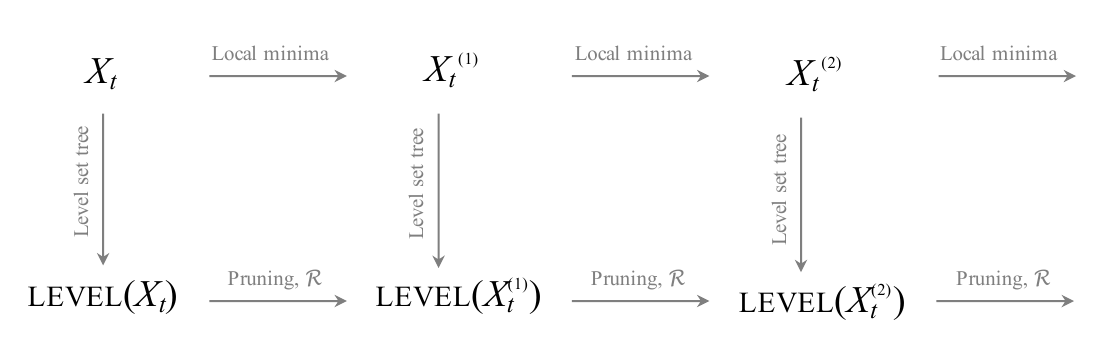}
\caption{Transition to the local minima of a function $X_t$ 
corresponds to pruning $\cR$ of the corresponding level set tree
$\textsc{level}(X_t)$.}
\label{fig:prune2}
\end{figure}

\begin{prop}[{{\bf Pruning for positive excursions, \cite{ZK12}}}]
\label{ts_prune}
Using the definitions of this section, the transition from a positive excursion 
$X_t$ to the respective excursion  $X^{(1)}_t$ 
of its local minima corresponds to pruning of the level-set tree $\textsc{level}(X_t)$.
This is illustrated in a diagram of Fig.~\ref{fig:prune2}.
Formally,
\[\textsc{level}\left(X^{(m)}_t\right) 
= \cR^m\left(\textsc{level}(X_t)\right), \forall m \ge 1.\]
\end{prop}

It is straightforward to formulate an analog of Prop.~\ref{ts_prune} 
without the excursion assumption (for continuous functions with a
finite number of local minima).
This, however, involves technicalities that are tangential to the essence
of this work and will be discussed elsewhere.

\subsection{Self-similarity for time series}\label{sec:level-set}
\begin{figure}[h] 
\centering\includegraphics[width=0.7\textwidth]{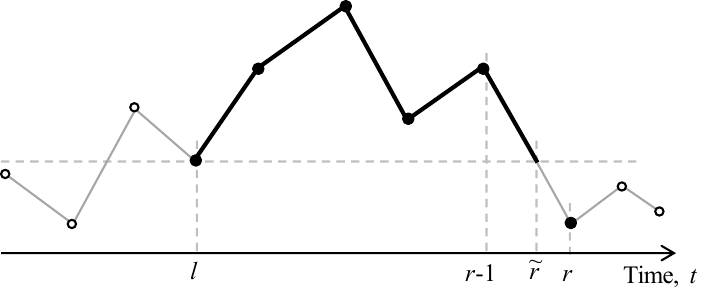}
\caption{Excursion of a time series: illustration.
The values of time series $X_k$ are shown by circles;
the circles that form the excursion on $[l,r]$ are filled.
The linear interpolation function $X_t$ is shown by solid line;
the excursion of $X_t$ on $[l,\tilde r]$ is shown in bold. 
}
\label{fig:ex}
\end{figure}

Consider a time series $X_k$, $k\in\mathbb{Z}$, with an atomless
distribution of values at each $k$.
Let $X_t$, $t\in\mathbb{R}$, be a continuous function of 
linearly interpolated values of $X_k$.
We define a {\it positive excursion} of $X_k$ as a fragment
of the time series on an interval $[l,r]$, $l,r\in\mathbb{Z}$
such that $X_l\ge X_r$ and $X_k>X_l$ for all $l<k<r$.
To each positive excursion of $X_k$ on $[l,r]$ corresponds a
positive excursion of $X_t$ on $[l,\tilde r]$, where $\tilde r\in(r-1,r]$ 
is such that $X_{\tilde r} = X_l$.
This construction is illustrated in Fig.~\ref{fig:ex}. 
The level set tree of a positive excursion of $X_k$ is that
of the corresponding positive excursion of $X_t$.

Proposition~\ref{ts_prune} and a comment after it suggest that the problem 
of finding self-similar trees with edge lengths is equivalent to finding 
{\it extreme-invariant processes}
\be
\label{d_inv}
X_k\stackrel{d}{=}\zeta\,X^{(1)}_k\quad {\rm for~some~} \zeta>0,
\ee
where $X_k$, $k\in\mathbb{Z}$, is a time series with an atomless value 
distribution at every $k$ and $X^{(1)}_k$ is the corresponding time series of 
local minima.
If $X_k$ satisfies \eqref{d_inv}, the level set tree of an excursion from $X_k$, considered
as an element of $\cL$, is self-similar according to Def.~\ref{def:distss}.
The next section describes a solution to \eqref{d_inv} that corresponds to $\zeta=2$.

\subsection{Self-similarity for random walks on $\mathbb{R}$}
\label{cgw}
Consider a random walk $\{X_k\}_{k\in\mathbb{Z}}$ with a homogeneous transition kernel $p(x,y)\equiv p(x-y)$, for any $x,y\in\mathbb{R}$, where $p(x)$ is an atomless density function.
A homogeneous random walk is called {\it symmetric} if $p(x)=p(-x)$ for all $x\in\mathbb{R}$.

\begin{lem}[{{\bf Pruning for random walks, \cite{ZK12}}}]
\label{inv_prune}
The following statements hold.
\begin{description}
  \item[a] The local minima of a homogeneous random walk $\{X_k\}_{k\in\mathbb{Z}}$ form a homogeneous random walk (with a different transition kernel in general). 
  \item[b] The local minima of a symmetric homogeneous random walk $\{X_k\}_{k\in\mathbb{Z}}$ form a symmetric homogeneous random walk (with a different transition kernel in general).
\end{description}
\end{lem}

The transition kernel of a symmetric random walk can be represented as the even part of 
a probability density function $f(x)$ with support in $\mathbb{R}_+$:
\[p(x)=\frac{f(x)+f(-x)}{2}.\] 
The following result describes the solution of the problem \eqref{d_inv}
in terms of the characteristic function of $f(x)$.

\begin{prop}[{{\bf Self-similarity for a symmetric homogeneous random walk, \cite{ZK12}}}]
\label{DSS}
The local minima of a symmetric homogeneous random walk  $\{X_k\}_{k\in\mathbb{Z}}$ 
with a transition kernel $~p(x)=\frac{f(x)+f(-x)}{2}~$ form a symmetric homogeneous random walk 
with a transition kernel
\[p^{(1)}(x)= \zeta^{-1}\,p(x/\zeta),\quad c>0\]
if and only if $\zeta=2$ and
\be
\label{laplace}
\Re\left[\widehat{f}(2s)\right]=\left|\frac{\widehat{f}(s)}{2-\widehat{f}(s)}\right|^2,
\ee
where $\widehat{f}(s)$ is the characteristic function of $f(x)$ and
$\Re[z]$ stays for the real part of $z\in\mathbb{C}$.
\end{prop}
A solution to \eqref{laplace} is given for example by an exponential density 
$f(x)=\phi_\lambda (x)$ of \eqref{exp} for any $\lambda>0$;
a detailed discussion of exponential kernels is given in Sect.~\ref{erw}.
A weaker, mean self-similarity of Defs.~\ref{ss1},~\ref{ss2} is satisfied in
any symmetric random walk, as discussed in the following statement.
  
\begin{thm}[{{\bf Mean self-similarity of a symmetric homogeneous random walk, \cite{ZK12}}}]
\label{T1}
The combinatorial level set tree $T=\textsc{shape}\left(\textsc{level}(X_t)\right)$ 
of a finite symmetric homogeneous random walk $X_k$ with $k=1,\dots,N$ is mean self-similar.
Specifically, for a uniform random complete subtree ${\sf subtree}_{K,{\sf k(T)}}\subset T$ 
of order $K<{\sf k}(T)$ the numbers $\tau_{i,j}^{(r)}$ of side-branches of order $i$
that merge the $r$-th branch of order $j$, with
$2\le j\le K$, in ${\sf subtree}_{K,{\sf k(T)}}$
are independent identically distributed random variables.
If $\tau_{i,j}$ is a random variable such that $\tau^{(r)}_{i,j}\stackrel{d}{=}\tau_{i,j}$, then
\be
\label{Tokunaga}
T_{i,j} := {\sf E}\left[\tau_{i,j}\right] = 2^{j-i-1}=:T_{j-i}.
\ee
Moreover, by the strong law of large numbers ${\sf k}(T)\overset{a.s.}{\to}\infty$  as $N \rightarrow \infty$,
and for any $1\le i<j$ we have
\[T_{i,j}\overset{\rm a.s.}{\longrightarrow} T_{j-i}=2^{j-i-1}
\quad{\rm as~}N\to\infty,\]
where $T_{i,j}$ can be computed over the entire $T$.
\end{thm}

\subsection{Exponential random walks}
\label{erw}
We call a homogeneous random walk {\it exponential} if its kernel is 
a mixture of exponential jumps constructed as follows
\[p(x)=\rho\,\phi_{\lambda_u}(x)+(1-\rho)\,\phi_{\lambda_d}(-x),
\quad 0\le \rho \le 1, \quad \lambda_u,\lambda_d>0,\]
where $\phi_{\lambda}$ is the exponential density with parameter $\lambda>0$,
\be
\label{exp}
\phi_\lambda (x)=
\begin{cases} 
\lambda e^{-\lambda x} &\text{ if } x \geq 0, \\ 
0 &\text{ if } x<0. 
\end{cases}
\ee
We refer to an exponential random walk by its parameter triplet $\{\rho,\lambda_u,\lambda_d\}$.
Each exponential random walk with parameters $\{\rho,\lambda_u,\lambda_d\}$ corresponds to a 
piece-wise linear function from $\cE(\mathbb{R})$ whose rises and falls have independent exponential lengths with parameters $(1-\rho)\lambda_u$ and $\rho\lambda_d$, respectively.
An exponential random walk is symmetric if and only if $\rho=1/2$ and $\lambda_u=\lambda_d$.

\begin{thm}[{{\bf Self-similarity of exponential random walks, \cite{ZK12}}}]
\label{eH}
Let $X_t$ be an exponential random walk with parameters $\{\rho,\lambda_u,\lambda_d\}$.
Then 
\begin{description}
\item[a]
The local minima of $X_t$ form a exponential random walk with parameters 
$\{\rho^*,\lambda_u^*,\lambda_d^*\}$ such that
\be 
\label{iteration}
\rho^*=\frac{\rho\,\lambda_d}{\rho\,\lambda_d+(1-\rho)\,\lambda_u}, \quad 
\lambda^*_d=\rho\lambda_d,~~\text{ and }~~
\lambda^*_u=(1-\rho)\lambda_u.
\ee
\item[b]
The exponential walk $X_t$ satisfies the self-similarity \eqref{d_inv}
if and only if it is symmetric, that is if $\rho=1/2$ and $\lambda_u=\lambda_d$.
\item[c]
The self-similarity \eqref{d_inv} is achieved after the first pruning, for the chain $X_t^{(1)}$ of the local minima,
if and only if the walk's increments have zero mean, $\rho\,\lambda_d = (1-\rho)\,\lambda_u$. 
\end{description}
\end{thm}

Recall that $\mathcal{GW}(p_0,p_2)$ denotes the space of binary Galton-Watson trees 
with termination probability $p_0$ and split probability $p_2$ (see Example~\ref{ex1}).

\begin{Def} [{{\bf Exponential binary Galton-Watson tree, \cite{Pitman}}}]
\label{def:binary}
We say that a random {\it embedded binary tree} 
$T\in\cL_{\rm plane}$ is an exponential binary Galton-Watson tree ${\sf GW}(\lambda',\lambda)$,
for $0\le\lambda'<\lambda$, if  
\textsc{shape}($T$) is a binary Galton-Watson tree $\mathcal{GW}(p_0,p_2)$ with 
\[p_0=\frac{\lambda+\lambda'}{2\lambda},\quad
p_2 = \frac{\lambda-\lambda'}{2\lambda},\]
and given \textsc{shape}($T$), the edges of $T$ are sampled as independent 
exponential random variables with parameter $2\lambda$, i.e., with density $\phi_{2\lambda} (x)$.
\end{Def}

A connection between exponential random walks and Galton-Watson trees  
is given by the following well known result.
\begin{thm}{\rm \cite[Lemma 7.3]{Pitman},\cite{LeGall93,NP}}
\label{Pit7_3}
Consider a random excursion $Y_t$ in $\cE^{\rm ex}$.
The level set tree $\textsc{level}(Y_t)$ 
is an exponential binary Galton-Watson tree ${\sf GW}(\lambda',\lambda)$ 
if and only if the rises and falls of $Y_t$, excluding the last fall, are 
distributed as independent exponential random variables with parameters $(\lambda+\lambda')$
and $(\lambda-\lambda')$, respectively, for some $0\le \lambda' < \lambda$.
Equivalently, the level set tree of a homogeneous random walk
is a binary Galton-Watson tree ${\sf GW}(\lambda',\lambda)$ if
and only if $Y_t$, as an element of $\cE^{\rm ex}$,  corresponds to an excursion of an exponential walk with 
parameters $\{\rho,\lambda_u,\lambda_p\}$
such that
$(1-\rho)\lambda_u = \lambda+\lambda'$ and
$\rho\lambda_d = \lambda-\lambda'.$
\end{thm}

\noindent We emphasize the following direct consequence of Thms.~\ref{eH}(a) and \ref{Pit7_3}.
\begin{cor}
\label{cor:GW}
Consider a critical binary Galton-Watson tree with independent exponential
lengths, $T = {\sf GW}(0,\gamma)$. The following statements hold:
\begin{description}
\item[a] The Harris path of $\cR^{k}(T)$ for any $0\le k<{\sf k}(T)$
corresponds to a positive excursion of a symmetric exponential random walk with parameters 
$\left\{\frac{1}{2},2^{1-k}\gamma,2^{1-k}\gamma\right\}$, or, equivalently, 
$\cR^{k}(T) = {\sf GW}\left(0,2^{-k}\gamma\right)$.
\item[b] The length of any branch of order $j\ge 1$ in $T$ has exponential distribution
with parameter $2^{2-j}\,\gamma$. The lengths of branches (of all orders) are independent.
\end{description}
\end{cor}

\section{Hierarchical branching process}
\label{HBP}
The results of previous section concern a very narrow class of mean self-similar
trees -- those with $T_k=2^{k-1}$.
Among such trees, the self-similarity is established 
only for the critical binary Galton-Watson tree ${\sf GW}(0,\gamma)$ 
with independent exponential edge lengths, i.e., continuous parameter Galton-Watson 
binary branching Markov processes; 
this case corresponds to the scaling exponent $\zeta=2$. 
Here we construct a branching process that generates self-similar 
trees for an arbitrary sequence $T_k\ge0$ and for any $\zeta > 0$; 
it includes the critical binary Galton-Watson tree as a special case.

\subsection{Definition and main properties}
\label{HBP:def}
Consider a probability mass function $\{p_K\}_{K\ge1}$, 
a sequence $\{T_k\}_{k\ge 1}$ of nonnegative Tokunaga coefficients, 
and a sequence $\{\lambda_j\}_{j\ge 1}$ of positive termination rates.
A multi-type branching process $S(t)$ starts with a root branch of Horton-Strahler
order $K\ge 1$ with probability $p_K$. 
Every branch of order $j\le K$ produces offspring of order $i<j$ with rate $\lambda_j T_{j-i}$. 
A branch of order $j$ terminates with rate $\lambda_j$.
After termination, a branch 
of order $j\geq 2$ splits into two branches of order $j-1$. 
A branch of order $j=1$ terminates without leaving offspring.  
The branching history of $S(t)$ creates a random binary tree $T[S]$ in the space $\cL$ of binary trees with edge lengths and no planar embedding.
The process is uniquely specified by the triplet
\[S(t) =\left(\{T_k\},\{\lambda_j\},\{p_K\}\right).\]

\begin{prop}[{{\bf Side-branching in hierarchical branching process}}]
\label{HBP:branch}
Consider a hierarchical branching process $S(t)=\left(\{T_k\},\{\lambda_j\},\{p_K\}\right)$.
For any branch $b\subset T[S]$ of order $K+1\ge 2$, let $m_i:=m_i(b)\ge 0$ be the
number of its side branches of order $i=1,\dots,K$, and $m:=m(b)=m_1+\dots+m_{K}$ be
the total number of the side branches.
Let $l_i:=l_i(b)$ be the lengths of $m+1$ edges within $b$, counted sequentially from 
the initial vertex,  
and $l:=l(b)=l_1+\dots+l_{m+1}$ be the total branch length.
Then the following statements hold:
\begin{enumerate}

\item The total numbers $m(b)$ of side branches within different branches of order $K+1$
are i.i.d. random variables with a common geometric distribution:
\be
\label{mdist}
{\sf P}\big(m=\kappa\big) = q(1-q)^\kappa \quad\text{with}\quad q=\frac{1}{1+T_1+\dots+T_{K}}, \quad \kappa=0,1,\hdots.
\ee

\item The number $m_i$ of side branches of order $i$ has geometric distribution:
\be
\label{midist}
{\sf P}\big(m_i=\kappa\big) = q_i(1-q_i)^\kappa \quad\text{with}\quad q_i=\frac{1}{1+T_{K-i+1}}, \quad \kappa=0,1,\hdots.
\ee

\item Conditioned on the total number $m$ of side branches, the distribution 
of $\{m_i\}$ is multinomial with $m$ trials and success probabilities
\be
\label{multinomial}
{\sf P}(\text{side~branch~has~order~}i)=\frac{T_{K-i+1}}{T_1+\dots+T_K}.
\ee
The side branch order vector ${\sf k}=({\sf k}_1,\dots,{\sf k}_m)$, where the side branches
are labeled sequentially starting from the initial vertex of $b$, is obtained
from the sequence
\[{\sf orders}=\underbrace{(1,\dots,1}_{m_1\text{~times}},\underbrace{2,\dots,2}_{m_2\text{~times}},
\dots\underbrace{K,\dots,K)}_{m_K\text{~times}}\]
by a uniform random permutation $\sigma_m$ of indices $\{1,\dots,m\}$:
\[{\sf k} = {\sf orders}\circ\sigma_m.\]

\item The branch length $l$ has exponential
distribution with rate $\lambda_{K+1}$, independent of the lengths of
any other branch (of any order). 
The corresponding edge lengths $l_i$ are i.i.d. random variables; they have
a common exponential distribution with rate
\be
\label{edges}
\lambda_{K+1}(1+T_1+\dots+T_K).
\ee
\end{enumerate}
\end{prop}
\noindent
\begin{proof}
All the properties readily follow from process construction.
\end{proof}

Proposition~\ref{HBP:branch} provides an alternative definition of the 
hierarchical branching process, and its construction -- via parts 
(1), (3), and (4) -- that does not require time-dependent simulations.

\begin{thm}[{{\bf Self-similarity of hierarchical branching process}}]
\label{main}
Consider a hierarchical branching process $S(t)=\left(\{T_k\},\{\lambda_j\},\{p_K\}\right)$
and let $T:=T[S]$ be the tree generated by $S(t)$.
The following statements hold.
\begin{enumerate}

\item The combinatorial part of $T$ is mean self-similar
(according to Def.~\ref{ss1},\ref{ss2}) with Tokunaga coefficients $\{T_k\}$.

\item The combinatorial part of $T$ is self-similar
(according to Def.~\ref{def:ss}) with Tokunaga coefficients $\{T_k\}$
if and only if 
\[p_K = p(1-p)^{K-1}\quad\]
for some $0<p< 1$.

\item The tree $T$ is self-similar (according to Def.~\ref{def:distss})
with scaling exponent $\zeta>0$ if and only if
\[p_K = p(1-p)^{K-1}\quad{\rm and}\quad\lambda_j = \gamma\,\zeta^{-j}\]
for some positive $\gamma$ and $0<p< 1$.
\end{enumerate} 
\end{thm}
\begin{proof}
By process construction, the tree $T$ is coordinated in shapes and lengths 
(according to Def.~\ref{def:distss}), with independent complete subtrees.

(1) Proposition~\ref{HBP:branch}, part (3) implies that the expected value of the number
$\tilde N_{i,j}$ of side branches of order $i\ge 1$ within a branch of order $j>i$ is 
given by ${\sf E}\left(\tilde N_{i,j}\right)=T_{j-i}$.
The mean self-similarity of Def.~\ref{ss1} with coefficients $T_k$ immediately follows, using 
a conditional argument as in \eqref{condarg}.

(2) Assume that $\textsc{shape}\left(T\right)$ is self-similar. 
A geometric distribution of orders is then established in Prop.~\ref{geom:p}.
Inversely, a geometric distribution of orders ensures that the total mass
$\mu\left(\cH_K\right)$, $K\ge 1$, is invariant with respect to pruning.
The conditional distribution of trees of a given order is completely specified
by the side branch distribution, described in Proposition~\ref{HBP:branch},
parts (1)-(3). 
Consider a branch of order $K+1$, $K\ge 1$. 
Pruning decreases the orders of this branch, and all its side branches, by unity.
Pruning eliminates a random geometric number $m_1$ of side-branches of order 1
from the branch.
It acts as a thinning (with removal probability $T_{K}/(T_1+\dots+T_K)$)
on the total side branch count $m$.
Accordingly, the total side branch count after pruning has geometric distribution
with success probability
\[q^{\cR} = \frac{1}{1+T_1+\dots T_{K-1}}.\]
The order assignment among the remaining side branches of orders $i=1,\dots,K-1$ 
is done according to multinomial distribution with probabilities proportional to $T_{K-i}$.
This coincides with the side branch structure in the original tree, hence completing 
the proof of (2).

(3) Having proven (2), it remains to prove the statement for the length structure of the tree.
Assume that $T$ is self-similar
with scaling exponent $\zeta$. 
The branches of order $j\ge 2$ become branches of order $j-1$ after
pruning, which necessitates $\lambda_{j} = \zeta\,\lambda_{j-1}$.
Inversely, pruning acts as a thinning on the side branches within a
branch of order $K+1$, eliminating the side branches of order ${\sf k}=1$. 
Accordingly, the spacings between the remaining side branches are
exponentially distributed with a decreased rate
\[\lambda_{K+1}(1+T_1+\dots+T_{K-1})=\zeta\,\lambda_K(1+T_1+\dots+T_{K-1}).\]
Comparing this with \eqref{edges}, and recalling the self-similarity of $\textsc{shape}\left(T\right)$, we conclude that Def.~\ref{def:distss}
is satisfied with scaling exponent $\zeta$.
\end{proof}

\subsection{Hydrodynamic limit}
\label{HBP:hydro}
Here we analyze the average numbers of branches of different orders in 
a hierarchical branching process, using a hydrodynamic limit.
Specifically, let $n\,x^{(n)}_j(s)$ be the number of branches of order $j$ at time 
$s$ observed in $n$ independent copies of the process $S$. 
Let $N_j(s)$ be the number of branches of order $j\ge 1$ in the process $S$
at instant $s\ge 0$. 
We observe that, by the law of large numbers,
\[ x^{(n)}_j(s)\stackrel{\text{a.s.}}{\longrightarrow}{\sf E}\left(N_j(s)\right)= :x_j(s).\]

\begin{thm}[{{\bf Hydrodynamic limit for branch dynamics}}]
\label{thm:hydro}
Suppose that the following conditions are satisfied:
\be\label{eq:L}
L:=\limsup_{k\to\infty} T_k^{1/k}<\infty,
\ee
and 
\be\label{con}
\sup\limits_{j \geq 1} \lambda_j <\infty,\quad
\limsup\limits_{j\to\infty} \lambda_j ^{1/j}\le 1/L.
\ee
Then, for any given $T>0$, the empirical process 
\[x^{(n)}(s)=\Big(x^{(n)}_1(s),x^{(n)}_2(s),\hdots \Big)^T, \qquad s \in [0,T],\] 
converges almost surely, as $n\to\infty$, to the process 
\[x(s)=\Big(x_1(s),x_2(s),\hdots \Big)^T, \qquad s \in [0,T],\] 
that satisfies
\begin{equation}\label{GODE}
\dot{x}=\mathbb{G}\Lambda x \quad
\text{ with the initial conditions } \quad x(0)=\pi:=\sum\limits_{K=1}^\infty p_K e_K,
\end{equation}
where $\Lambda=\text{diag}\{\lambda_1,\lambda_2,\hdots\}$ is a diagonal operator with 
the entries $\lambda_1,\lambda_2,\hdots~$,
$e_i$ are the standard basis vectors, and
\begin{equation} \label{gen}
\mathbb{G}:=\left[\begin{array}{ccccc}-1 & T_1+2 & T_2 & T_3 & \hdots  \\0 & -1 & T_1+2 & T_2 & \hdots  \\0 & 0 & -1 & T_1+2 & \ddots  \\0 & 0 & 0 & -1  & \ddots \\\vdots & \vdots & \ddots & \ddots  & \ddots \end{array}\right].
\end{equation}
\end{thm}
\begin{proof} 
The process $~x^{(n)}(s)~$ evolves according to the transition rates 
$$q^{(n)}(x,x+\ell)=n\beta_\ell\left({1 \over n} x\right)$$ 
with
$$\beta_\ell(x)=\begin{cases}
      \lambda_1x_1 & \text{  if } \ell=-e_1, \\
      \lambda_{i+1}x_{i+1} & \text{  if } \ell=2e_i-e_{i+1}, i\ge 1,\\
      \sum\limits_{j=i+1}^\infty \lambda_j T_{j-i}x_j & \text{ if } \ell=e_i, i\ge 1.
\end{cases}$$
Here the first term reflects termination of branches of order 1;
the second term reflects termination of branches of orders $i+1>1$, each of 
which results in creation of two branches of order $i$; and
the last term reflects side-branching.
Thus, the infinitesimal generator of the stochastic process $x^{(n)}(s)$ is
\begin{eqnarray}\label{eq:Ln}
L_nf(x) & = & n\lambda_1 x_1\left[f\left(x-{1 \over n}e_1\right)-f(x)\right]  +\sum\limits_{i=1}^\infty n\lambda_{i+1}x_{i+1}\left[f\left(x-{1 \over n}e_{i+1}+{2 \over n}e_i\right)-f(x)\right]  \nonumber \\
&  &  ~~+  \sum\limits_{i=1}^\infty\left(\sum\limits_{j=i+1}^\infty n \lambda_jT_{j-i}x_j\right) \left[f\left(x+{1 \over n}e_i\right)-f(x)\right].
\end{eqnarray}
Let
$$F(x):=\sum\limits_\ell \beta_\ell(x)=-\lambda_1 x_1e_1 +\sum\limits_{i=1}^\infty \lambda_{i+1}x_{i+1}(2e_i-e_{i+1})+\sum\limits_{i=1}^\infty \left(\sum\limits_{j=i+1}^\infty \lambda_j T_{j-i}x_j\right) e_i.$$
The  convergence result of Kurtz (\cite[Theorem 2.1, Chapter 11]{EK86}, \cite[Theorem 8.1]{Kurtz81}) extends (without changing the proof) to the Banach space $\ell^1(\mathbb{R})$ provided the same conditions are satisfied 
for $\ell^1(\mathbb{R})$ as for $\mathbb{R}^d$ in the theorem of Kurtz. 
Specifically, we require that for a compact set $\mathcal{C}$ in $\ell^1(\mathbb{R})$, 
\begin{equation}\label{eq:beta}
\sum_\ell \|\ell\|_1 \sup_{x \in \mathcal{C}} \beta_\ell(x) <\infty,
\end{equation} 
and there exists $M_\mathcal{C} >0$ such that 
\begin{equation}\label{eq:Lipschitz}
\|F(x)-F(y)\|_1 \leq M_\mathcal{C} \|x-y\|_1, \qquad x,y \in \mathcal{C} .
\end{equation}
Here the condition (\ref{eq:beta}) follows from
$$\sum_i \sup_{x \in \mathcal{C}} |\lambda_i x_i| <\infty \qquad \text{ and } \qquad \sum_i \sup_{x \in \mathcal{C}} \sum\limits_{j=i+1}^\infty \lambda_j T_{j-i}|x_j| <\infty,$$
which in turn follow from conditions \eqref{con}. 
Similarly, Lipschitz conditions (\ref{eq:Lipschitz}) are satisfied in 
$\mathcal{C}$ due to conditions \eqref{con}.
Thus, by Kurtz (\cite[Theorem 2.1, Chapter 11]{EK86}, \cite[Theorem 8.1]{Kurtz81}),  
the process $x^{(n)}(s)$ converges almost surely to $x(s)$ that satisfies $\dot{x}=F(x)$,
which expands as the following system of ordinary differential equations:
\begin{equation}\label{gFODEsys}
\begin{cases}
x'_1(s) & =  -\lambda_1 x_1 +\lambda_2 (T_1+2)x_2+\lambda_3 T_2x_3+\hdots \\ 
x'_2(s) & =  -\lambda_2 x_2 +\lambda_3 (T_1+2)x_3+\lambda_4 T_2x_4+\hdots \\  
& \vdots \\
x'_k(s) & =  -\lambda_k x_k +\lambda_{k+1} (T_1+2)x_{k+1}+\lambda_{k+2} T_2x_{k+2}+\hdots \\ 
& \vdots   
\end{cases}
\end{equation}
with the initial conditions $x(0)=\lim\limits_{n \rightarrow \infty} x^{(n)}(0)=\pi:=\sum\limits_{K=1}^\infty p_K e_K$ by the law of large numbers. 
Finally, we observe that $\|\pi\|_1=1$, and conditions \eqref{con} imply that 
$\mathbb{G}\Lambda$ is a bounded operator in $\ell^1(\mathbb{R})$.
\end{proof}

\subsection{Criticality and time invariance}
\label{sec:width}
Assume that the hydrodynamic limit $x(s)$, and hence the averages $x_j(s)$, exist. Let $\pi=\sum\limits_{K=1}^\infty p_K e_K$.
Then one can consider the average progeny of the process, that is the average 
number of branches of any order alive at instant $s$:
$$C(s)=\sum\limits_{j=1}^\infty x_j(s)=\Big\|e^{\mathbb{G}\Lambda s} \pi \Big\|_1.$$
In hydrological literature, an empirical version of the process $C(s)$ is 
called the {\it width function} of a tree $T[S]$.

\begin{Def}\label{criticalpi}
A hierarchical branching process $S(s)$ is said to be 
{\it critical} if and only if the width function $C(s)=1$ for all $s \geq 0$.
\end{Def}

\begin{Def}
A hierarchical branching process $S(s)$ is said to be
{\it time-invariant} if and only if 
\begin{equation}\label{eq:fipc}
e^{\mathbb{G}\Lambda s}\pi=\pi\quad\text{ for all}\quad s\ge 0.
\end{equation}
\end{Def}

\begin{prop}\label{ti2c}
Suppose that the hydrodynamic limit $x(s)$ exists, and $\pi$ is time-invariant. 
Then the process $S(s)$ is critical.
\end{prop}
\begin{proof}
$C(s)=\|x(s)\|_1=\| e^{\mathbb{G}\Lambda s} \pi \|_1=\| \pi \|_1=1.$
\end{proof}

Let $~\hat{t}(z) = -1+2z+\sum_j z^j\,T_j~$ for $|z| <1/L$, where $L$ is defined in (\ref{eq:L}).
Observe that there is a unique real root $w_0$ of $\hat{t}(z)$ within $(0,\frac{1}{2}]$.
We formulate our results in terms of the {\it Horton exponent} $R:=w_0^{-1}$ (e.g., \cite{Pec95,KZ15t}). 

\begin{prop}
\label{critcond}
Suppose 
$\Lambda\,\pi$ is a constant multiple of the geometric vector 
$v_0=\sum\limits_{K=1}^\infty R^{-K} e_K$. 
Then the process $S(s)$ is time-invariant.
\end{prop}
\begin{proof}
Observe that since $\hat{t}\left(R^{-1}\right)=0$ and $\mathbb{G}$ is a Toeplitz operator,
$$\mathbb{G}v=\hat{t}(w)v \quad \text{ for }~v=\sum\limits_{K=1}^\infty w^K e_K, ~~|w|<L.$$
and 
$$\mathbb{G}v_0=\hat{t}\left(R^{-1}\right)v_0=0 \quad \text{ for }~v_0:=
\sum\limits_{K=1}^\infty R^{-K} e_K.$$
Hence
$\mathbb{G}\Lambda \pi = \hat{t}\left(R^{-1}\right) \Lambda \pi=0$
and
\[e^{\mathbb{G}\Lambda s}\pi = \pi +\sum_{m=1}^{\infty}\frac{s^m}{m!} (\mathbb{G}\Lambda)^m\pi=\pi.\]
\end{proof}

\begin{Rem}
Proposition~\ref{critcond} states that the condition
\be\label{lpR}
\lambda_K\,p_K = b\,R^{-K}, K\ge 1
\ee
is sufficient for time-invariance, for any proportionality constant $b>0$.
This implies that a time-invariant process can be constructed for 
\begin{itemize}
\item[(i)] an arbitrary sequence of Tokunaga coefficients $\{T_k\}$
satisfying \eqref{eq:L} -- by selecting $\lambda_K\,p_K = b\,R^{-K}$;
\item[(ii)] arbitrary sequences $\{T_k\}$ satisfying \eqref{eq:L} and $\{p_K\}$ -- by selecting
$\lambda_K = b\,R^{-K}\,p_K^{-1}$;
\item[(iii)] arbitrary sequences $\{T_k\}$ satisfying \eqref{eq:L} and $\{\lambda_K\}$ -- by selecting
$p_K = b\,R^{-K}\,\lambda_K^{-1}$.
\end{itemize}
At the same time, arbitrary sequences $\{\lambda_K\}, \{p_K\}$ 
will not, in general, satisfy \eqref{lpR} and hence will not
correspond to a time-invariant process.  
\end{Rem}

\subsection{Criticality and time-invariance in a self-similar process}
\label{HBP:ss}
A convenient characterization of criticality can be established for 
self-similar hierarchical branching processes.
Recall that by Theorem~\ref{main}, part (3), a self-similar process
$S(s)$ is specified by parameters $\gamma>0$, $0<p<1$ and length self-similarity 
constant $\zeta>0$ such that $p_K=p(1-p)^{K-1}$ and $\lambda_j = \gamma\,\zeta^{-j}$. 
We refer to a self-similar process by its parameter triplet,
$S(s)=S_{p,\gamma,\zeta}(s)$, and denote the respective width function 
by $C_{p,\gamma,\zeta}(s)$.
Observe that in the self-similar case the first of the conditions~\eqref{con}
is equivalent to $\zeta\ge 1$, and the second is equivalent to $\zeta\ge L$.
Hence, the conditions \eqref{con} are equivalent to
$\zeta \ge  1\vee L$.

\begin{thm}[{{\bf Width function of a self-similar process}}]
\label{wfss}
Consider a self-similar process $S_{p,\gamma,\zeta}(s)$
with $0<p<1$, $\gamma>0$.
Suppose that \eqref{eq:L} is satisfied and $\zeta \ge 1\vee L$.
Then
$$C_{p,\gamma,\zeta}(s)~\begin{cases}
   \text{decreases}    & \text{ if } p>1-{\zeta \over R}, \\
   = 1 & \text{ if  } p=1-{\zeta \over R}, \\
   \text{increases} & \text{ if  } p<1-{\zeta \over R}.
\end{cases}$$
\end{thm}
\begin{proof}
The choice of the limits for $\zeta$ ensures that the conditions \eqref{con}
are satisfied and hence, by Theorem~\ref{thm:hydro}, the hydrodynamic limit $x(s)$
exists
and  the width function $C_{p,\gamma,\zeta}(s)$ is well defined.
Now we have
\[~\Lambda \pi={\gamma p \over 1-p}\sum\limits_{K=1}^\infty \big(\zeta^{-1}(1-p)\big)^K e_K,\]
and  therefore 
\begin{equation}\label{eq:GLt}
\mathbb{G}\Lambda \pi = \hat{t}\big(\zeta^{-1} (1-p)\big) \Lambda \pi.
\end{equation} 
Iterating recursively, we obtain
$$(\mathbb{G}\Lambda)^2 \pi = \hat{t}\big(\zeta^{-1} (1-p)\big) \mathbb{G}\Lambda^2\pi= \hat{t}\big(\zeta^{-1} (1-p)\big) \hat{t}\big(\zeta^{-2} (1-p)\big) \Lambda^2 \pi,$$
and in general,
$$(\mathbb{G}\Lambda)^m \pi = \hat{t}\big(\zeta^{-1}(1-p)\big) \mathbb{G}\Lambda^m\pi= \left[\prod\limits_{i=1}^m \hat{t}\big(\zeta^{-i} (1-p)\big)\right] \Lambda^m \pi.$$
Thus, taking $x(0)=\pi$,
\begin{equation} \label{eq:solx}
x(s)=e^{\mathbb{G}\Lambda s} \pi=\pi+\sum\limits_{m=1}^\infty {s^m \over m!} \left[\prod\limits_{i=1}^m \hat{t}\big(\zeta^{-i}(1-p)\big)\right] \Lambda^m \pi.
\end{equation}
The width function for the given values of $p \in (0,1)$, $\gamma>0$ and $\zeta \ge 1$ can therefore be expressed as
\begin{align} \label{eq:Cpgz}
C_{p,\gamma,\zeta}(s)=\sum\limits_{j=1}^\infty x_j(s) &=1+\sum\limits_{m=1}^\infty {s^m \over m!} \left[\prod\limits_{i=1}^m \hat{t}\big(\zeta^{-i}(1-p)\big)\right] \sum\limits_{j=1}^\infty \big(\Lambda^m \pi\big)_j \nonumber \\
&=1+\sum\limits_{m=1}^\infty {\big(s\gamma/\zeta \big)^m \over m!} \left[\prod\limits_{i=1}^m \hat{t}\big(\zeta^{-i}(1-p)\big)\right] {p \over 1-\zeta^{-m}(1-p)} 
\end{align}
as $~\sum\limits_{j=1}^\infty \big(\Lambda^m \pi\big)_j=\sum\limits_{j=1}^\infty \lambda_j^m \pi_j=\sum\limits_{j=1}^\infty \gamma^m \zeta^{-jm} p(1-p)^{j-1}=\gamma^m \zeta^{-m} {p \over 1-\zeta^{-m}(1-p)}$.

Next, notice that by letting $p'=1-\zeta^{-1}(1-p)$, we have from (\ref{eq:Cpgz}) and the uniform convergence of the corresponding series for any fixed $M>0$ and $s \in [0,M]$, that
\begin{equation}\label{eq:Cprime}
{d \over ds}C_{p,\gamma,\zeta}(s)={\gamma \over \zeta} \hat{t}(1-p') C_{p',\gamma,\zeta}(s) \quad \text{ with }~ C_{p,\gamma,\zeta}(0)=C_{p',\gamma,\zeta}(0)=1.
\end{equation}
Observe that  $\zeta\geq 1$ implies $p' \geq p$ and $~C_{p',\gamma,\zeta}(s)\leq C_{p,\gamma,\zeta}(s)$. Also, observe that  
$$\hat{t}(1-p') ~\begin{cases}
   <0   & \text{ if } p>1-{\zeta \over R} \\
   =0 & \text{ if  } p=1-{\zeta \over R} \\
   >0 & \text{ if  } p<1-{\zeta \over R} 
\end{cases}$$
as $\hat{t}$ is an increasing function on $[0,\infty)$ and  $\hat{t}\big(1/R\big)=0$. 
This leads to the statement of the theorem.
\end{proof}

\begin{Rem}\label{spec_case}
If $\zeta=1$, equation (\ref{eq:Cprime}) 
implies $C_{p,\gamma,1}(s)=\exp\big\{s\gamma \hat{t}(1-p)\big\}$
and hence
$$C_{p,\gamma,1}(s)~\begin{cases}
   \text{ exponentially decreases}    & \text{ if } p>1-R^{-1}, \\
   =1 \text{ for all }s \geq 0 & \text{ if  } p=1-R^{-1}, \\
   \text{ exponentially increases} & \text{ if  } p<1-R^{-1}.
\end{cases}$$
This case is further examined in Sect.~\ref{DTM}.
In general, the width function $C_{p,\gamma,\zeta}(s)$ may increase 
sub-exponentially for $p<1-{\zeta \over R}$. 
For example, if there is a nonnegative integer $d$ such that $\zeta^{d+1} < R$,
then for $p=1-{\zeta^{d+1} \over R}$ we have
$\hat{t}\big(\zeta^{-d-1} (1-p)\big)=0$.
Hence, \eqref{eq:solx} implies that $C_{p,\gamma,\zeta}(s)$ is a polynomial of degree $d$.
\end{Rem}

\begin{thm}[{{\bf Criticality of a self-similar process}}]
\label{main2}
Consider a self-similar process $S_{p,\gamma,\zeta}(s)$
with $0<p<1$, $\gamma>0$.
Suppose that \eqref{eq:L} is satisfied and $\zeta \ge 1\vee L$. 
Then the following conditions are equivalent:
\begin{itemize}
\item[(i)]
The process is {\it critical}.
\item[(ii)]
The process is {\it time-invariant}.
\item[(iii)]
The following relations hold:
$\zeta<R\quad\text{and}\quad p=p_c: = 1-\frac{\zeta}{R}.$
\end{itemize}
\end{thm}
\begin{proof}
(i)$\leftrightarrow$(iii) is established in Theorem~\ref{wfss}.
(ii)$\rightarrow$(i) is established in Prop~\ref{ti2c}.
(iii)$\rightarrow$(ii):
Observe that $\hat{t}\left(\zeta^{-1}(1-p)\right)=\hat{t}\left(R^{-1}\right)=0$.
Time invariance now follows from (\ref{eq:solx}).
\end{proof}

\begin{Rem}
In a self-similar process the sequences $\lambda_K$ and $p_K$ are geometric
such that (Thm.~\ref{main}) 
\[\lambda_K\,p_K = \frac{\gamma\,p}{1-p}\left(\frac{1-p}{\zeta}\right)^K\]
for some $0<p<1$, $\gamma>0$, and $\zeta\ge 1\vee L$.
Hence, a time-invariant process can be constructed, according to
Prop.~\ref{critcond} and \eqref{lpR}, by selecting
any sequence $\{T_k\}$ that corresponds to
\[R = \zeta\,(1-p)^{-1}.\]
Theorem~\ref{main2} states that this is the only possible way
to construct a time-invariant process, given that the process
is self-similar.
\end{Rem}

\subsection{A closed form solution for the case of equally distributed branch lengths}
\label{DTM}
Observe that if $\Lambda=\gamma I$, then
$~x(s)=e^{s\gamma\hat{t}(1-p)}\pi~$ and $~C(s)=\|x(s) \|_1=e^{s\gamma \hat{t}(1-p)s}$. 

Consider a hierarchical branching process with $\Lambda=I$ and $x(0)=e_K$ for a given integer $K \geq 1$. Here the system of equation (\ref{gFODEsys}) is finite dimensional,
\begin{equation}\label{FODEsys}
\begin{cases}
x'_1(s) & =  -x_1 +(T_1+2)x_2+T_2x_3+\hdots +T_{K-1}x_K \\ 
x'_2(s) & =  -x_2 +(T_1+2)x_3+T_2x_4+\hdots +T_{K-2}x_K \\ 
& \vdots \\
x'_{K-1}(s) & =  -x_{K-1} +(T_1+2)x_K \\
x'_K(s) & =  -x_K 
\end{cases}
\end{equation}
with the initial conditions $x(0)=e_K$. 

\bigskip
\noindent
Define a sequence $t(j)$ as
$$t(0)=-1,~t(1)=T_1+2,~\text{ and }t(j)=T_j \text{ for }j \geq 2,$$
and let $y(s)=e^s x(s)$. Then (\ref{FODEsys}) becomes
\begin{equation}\label{yODEsys}
\begin{cases}
y'_1(s) & =  t(1)y_2+t(2)y_3+\hdots +t(K-1)y_K  \\ 
y'_2(s) & =  t(1)y_3+t(2)y_4+\hdots +t(K-2)y_K  \\ 
& \vdots \\
y'_{K-2}(s) & =  t(1)y_{K-1}+t(2)y_K \\
y'_{K-1}(s) & =  t(1)y_K \\
y'_K(s) & =  0 
\end{cases}
\end{equation}
with the initial conditions $y(0)=e_K$. The ODEs (\ref{yODEsys}) can be solved recursively in a reversed order of equations in the system obtaining for for $m=1,\hdots,K-1$,
$$y_{K-m}(s)=\sum\limits_{n=1}^m \left(\sum\limits_{\substack{i_1,\hdots,i_n \geq 1\\ i_1+\hdots+i_n=m}} t(i_1)\cdot \hdots \cdot t(i_n) \right){s^n \over n!}.$$
Let $\delta_0(j)=I_{\{j=0\}}$ be the Kronecker delta function. 
Then we arrive with the closed form solution
\begin{align}\label{eq:convt}
x_{K-m}(s)=e^{-s}y_{K-m}(s) & =e^{-s}\sum\limits_{n=1}^{\infty} \underbrace{(t+\delta_0)*(t+\delta_0)*\hdots*(t+\delta_0)}_{n \text{ times}} (m) {s^n \over n!}.
\end{align} 
Observe that if we randomize the orders of trees by assigning an order $K$ to a tree with geometric probability $p_K=p(1-p)^{K-1}$, 
then the above closed form expression (\ref{eq:convt}) would yield an expression 
for the width function that was observed in Remark \ref{spec_case} of this section:
\begin{eqnarray*}
C(s) & = & e^{-s}+e^{-s}\sum\limits_{n=1}^\infty \sum_{m=1}^\infty (1-p)^m ~\underbrace{(t+\delta_0)*(t+\delta_0)*\hdots*(t+\delta_0)}_{n \text{ times}} (m) ~{s^n \over n!} \\ 
&=&  e^{-s}+e^{-s}\sum\limits_{n=1}^\infty \Big(\hat{t}(1-p)+1\Big)^n ~{s^n \over n!}   
=  \exp\left\{s\hat{t}(1-p)\right\}.
\end{eqnarray*}

\subsection{Critical Galton-Watson process}
\label{sec:cGW} 
The critical binary Galton-Watson process plays an important role in 
theory and applications because of its multiple symmetries.
Burd, Waymire and Winn \cite{BWW00} have shown that the following three 
properties are equivalent for the binary Galton-Watson distributions $\mathcal{GW}(p_0,p_2)$:
(i) A distribution is prune-invariant;
(ii) A distribution is mean self-similar with ${\sf E}(T_{i,j}) = T_{j-i} = 2^{j-i-1};$
and (iii) A distribution is critical: $p_0=p_2=1/2$.
The Markov structure of the critical Galton-Watson trees ensures the existence 
of two other special properties:
(iv) Time-invariance (in discrete time): the forest of trees, obtained by removing 
the edges and the vertices below depth $d$, has the same frequency structure as the
original space $\mathcal{GW}(1/2,1/2)$;
and
(v) The forest of trees obtained by considering subtrees rooted
at every vertex of a random tree $T$ approximates the frequency 
structure of the entire space of trees when the order of $T$ increases.

The next results shows that the critical binary Galton-Watson tree is a 
special case of the hierarchical branching process.

\begin{thm}[{{\bf Critical Galton-Watson tree}}]
\label{main3}
A hierarchical branching process with parameters
\begin{equation}\label{eq:GW}
\lambda_j=\gamma 2^{2-j}, ~ p_K=2^{-K},~ \text{and}~ T_{k}=2^{k-1}
~\text{for~any~}\gamma>0
\end{equation} 
is distributionally equivalent to the critical binary Galton-Watson tree ${\sf GW}(0,\gamma)$
with i.i.d. edge lengths that have a common exponential distribution with rate $2\,\gamma$.
This is a self-similar, critical, and time-invariant process with
\[R = 4,\quad L = 2,\quad\text{and}\quad \zeta =2.\]
\end{thm}
 
\begin{proof}
Consider a tree $T={\sf GW}(0,\gamma)\in \cL$. 
By Corollary~\ref{cor:GW}, each branch of order $j$ in $T$ is exponentially 
distributed with parameter $\lambda_j=\gamma 2^{2-j}$, which matches the
branch length distribution in the hierarchical branching process \eqref{eq:GW}.
Furthermore, conditioned on $\cR^{i}(T)\ne \phi$ (which happens with a positive probability), 
we have $\cR^{i}(T)={\sf GW}(0,2^{-i}\gamma)$.
This means that the space $\cR^i(\cL)$ of pruned trees is a linearly scaled version of the 
original space $\cL$ (the same combinatorial structure, linearly scaled edge lengths).
Burd et al. \cite{BWW00} have shown that the total number of sub-branches 
within a branch of order $j\ge 2$ in $T$ is geometrically distributed over $\mathbb{Z}_+$
with mean $T_1+\dots+T_{j-1}=2^{j-1}-1$ (that is ${\sf P}(m) = 2^{1-j}\left(1-2^{1-j}\right)^m$
for $m\in \mathbb{Z}_+$),
where $T_{i,j}=T_{j-i}=2^{j-i-1}$.
The assignment of orders among the $m$ side-branches is done according to
the multinomial distribution with $m$ trials and success probabilities
$T_i/(T_1+\dots+T_{j-1})$, $i=1,\dots,j-1$.
This implies that, conditioning on a particular implementation of the pruned tree $\cR(T)$, 
the leaves of the original tree merge into every branch of the pruned tree as a Poisson 
point process with intensity $\gamma=\lambda_j T_{j-1}$.
Iterating this pruning argument, conditioning on the particular implementation of 
$\cR^i(T)={\sf GW}(0,2^{-i}\gamma)$, 
the branches of order $i$ merge into any branch of the pruned tree $\cR^i(T)$ as 
a Poisson point process with intensity $\gamma\,2^{1-i}=\lambda_j T_{j-i}$ for
every $j>i$.
Finally, the critical binary Galton-Watson space has $p_K=2^{-K}$ \cite{BWW00}. 
We, hence, conclude that a ${\sf GW}(0,\gamma)$ tree is distributionally identical 
to the hierarchical branching process with parameters (\ref{eq:GW}).

By Thm.~\ref{main}, the process~\eqref{eq:GW} is self-similar with
the scaling exponent $\zeta=2$.
Criticality and time-invariance follow from Prop.~\ref{main2}.
\end{proof}

\subsection{Critical Tokunaga processes}
\label{sec:Tok}
We introduce here a class of processes that extends 
the symmetries observed in the critical binary Galton-Watson 
tree with exponential edge lengths (where $\zeta=2$) to the general case of $\zeta\ge1$.
Specifically, consider a hierarchical branching process $S^{\rm Tok}(t;c,\gamma)$,
which we call the {\it critical Tokunaga branching process}, with parameters
\begin{equation}\label{eq:Tok}
\lambda_j=\gamma\,c^{2-j}, ~ p_K=2^{-K}, ~ \text{ and }~T_{k}=(c-1)\,c^{k-1}
~\text{for any }\gamma>0,~c\ge1.
\end{equation}

\begin{prop}
\label{Tok_tree}
The process $S^{\rm Tok}(t;c,\gamma)$ is a self-similar critical
time invariant process.
Independently of the process combinatorial shape, its edge lengths 
are i.i.d. exponential random variables with rate $\gamma c$.
In addition, we have
\[\hat{t}(z) = \frac{(1-2\,c\,z)(z-1)}{1-c\,z},~
R = w_0^{-1} = 2\,c,~\zeta = L = c,\text{ and }
p_c = 2^{-1}.\]
\end{prop}
\begin{proof}
Self-similarity follows from Thm.~\ref{main}.
Criticality and time-invariance are established in Prop.~\ref{main2}.
The edge lengths property follows from Prop.~\ref{HBP:branch}, part(4).
\end{proof}

\begin{Rem}
The condition $T_{i,i+k}=T_k=a\,c^{k-1}$ was first introduced in hydrology by
Eiji Tokunaga \cite{Tok78} in a study of river networks, hence the process 
name.
The additional constraint $a=c-1$ is necessitated here by the self-similarity
of tree lengths, which requires the sequence $\lambda_j$ to be geometric. 
The sequence of the Tokunaga coefficients then also has to be geometric, and
satisfy $a=c-1$, to ensure identical distribution of the edge lengths, 
see Prop.~\ref{HBP:branch}, part(4).
Interestingly, the constraint $a=c-1$ appears in the Random
Self-similar Network (RSN) model introduced by Veitzer
and Gupta \cite{VG00}, which uses a purely topological algorithm
of recursive local replacement of the network generators to generate
self-similar random trees.
The importance of the constraint $a=c-1$ in a combinatorial
situation is discussed in the next section. 
\end{Rem}

\section{Combinatorial critical Tokunaga process}
\label{sec:shape}
This section focuses on the combinatorial structure of a tree generated by the
critical Tokunaga process (leaving aside the edge lengths).
A general combinatorial reformulation of the hierarchical branching process 
can be found in \cite{KZ18}.

Consider a  combinatorially self-similar (according to Def. \ref{def:ss}) tree $T=\textsc{shape}\big(T[S]\big)\in\cT$ generated by
a hierarchical branching process $S$ with Tokunaga sequence $\{T_k\}$ 
and initial distribution $p_K=p(1-p)^{K-1}$.
Let random variable $K$ be the order of the tree $T$, and, conditioned on  $K>1$, let $K_a,K_b$ be the orders 
of its two subtrees, $T_a$ and $T_b$, rooted at the internal vertex closest to the root, randomly and uniformly permuted.
We call $T_a$ and $T_b$ the {\it principal} subtrees of $T$.
Observe that the pair $K_a,K_b$ uniquely defines the tree order $K$:
\[K=\begin{cases}
K_a\vee K_b,&\text{ if } K_a\ne K_b,\\
K_a+1,&\text{ if } K_a = K_b.
\end{cases}\]
Let $K_1\le K_2$ be the order statistics of $K_a,K_b$.
The joint distribution of $(K_1,K_2)$ is given by
\be
\label{joint}
{\sf P}\left(K_1=j,K_2=m|K=k\right)
=
\begin{cases}
\displaystyle\frac{1}{1+T_1+\dots+T_{k-1}}&\text{ if }j=m=k-1\\
\displaystyle\frac{T_{k-j}}{1+T_1+\dots+T_{k-1}}&\text{ if }j<m=k
\end{cases},
\ee
where
$${\sf P}(K=k|K>1)=(1-p)p^{k-2}.$$

\begin{prop}\label{subtrees}
Consider a critical Tokunaga process $S^{\rm Tok}(t;c,\gamma)$.
Then, conditioned on $K>1$, the marginal order distribution of $K_a$ coincides
with that of $K$: 
\be
\label{cT:marg}
{\sf P}(K_a=m ~|~K>1)=2^{-m} = p_m \quad \text{ for }\quad m\ge 1.
\ee
At the same time, the joint distribution of $(K_a,K_b)$ equals 
the product of the marginals,
\be
\label{cT:joint}
{\sf P}(K_a=m, ~K_b=j ~|~K>1)={\sf P}(K_a=m ~|~K>1){\sf P}(K_b=j ~|~K>1),
\ee
if and only if $c=2$.
\end{prop}

\begin{proof}
Observe that for $k>1$,
\begin{align*}
{\sf P}(K_a=m ~&|~K=k)\\
=&
\begin{cases} 
\frac{1}{2} \sum\limits_{j:j<k} {\sf P}(K_1=j,K_2=k|K=k)&\text{ if }m=k,\\
{\sf P}(K_1=K_2=k-1|K=k)+\frac{1}{2} {\sf P}(K_1=k-1,K_2=k|K=k)&\text{ if }m=k-1,\\
\frac{1}{2} {\sf P}(K_1=m,K_2=k|K=k)&\text{ if }m<k-1,
\end{cases}\\
=&
\left\{
\begin{array}{lll}
\displaystyle{1\over 2}{T_1+\dots+T_{k-1}  \over 1+T_1+\dots+T_{k-1}} &={1-c^{1-k}\over 2} &\text{ if }m=k,\\
\displaystyle\frac{1+\frac{1}{2}T_1}{1+T_1+\dots+T_{k-1}} &={c+1 \over 2}c^{1-k} &\text{ if }m=k-1,\\
\displaystyle{1\over 2}{T_{k-m}  \over 1+T_1+\dots+T_{k-1}} &={(c-1)c^{-m}\over 2} &\text{ if }m<k-1.
\end{array}
\right.
\end{align*}
Thus,
\begin{align*}
{\sf P}(K_a=m ~|~K>1)
=& (1-c^{1-m})2^{-m}+c^{-m}2^{-m}+{(c-1)c^{-m} \over 2} \sum\limits_{k=m+1}^\infty 2^{1-k}
= 2^{-m}.
\end{align*}
At the same time,
\begin{align*}
{\sf P}(K_a=m, ~K_b=j ~|~K>1)
=& 
\begin{cases}
(c-1)c^{-j}\,2^{-m} &\text{ if }j<m,\\
c^{-m}\,2^{-m}&\text{ if }j=m.
\end{cases}
\end{align*}
Hence, \eqref{cT:joint} holds if and only if $c=2$.
\end{proof}

\begin{Rem} 
Proposition~\ref{subtrees} asserts that the principal
subtrees in a random critical Tokunaga tree are dependent, except 
the critical Galton-Watson case.
This implies that, in general, non-overlapping subtrees within
a critical Tokunaga tree are dependent.
Accordingly, the increments of the Harris path $H$ of a critical Tokunaga 
process have (long-range) dependence.
The only exception is the case $c=2$ that was discussed in Sect.~\ref{erw}.
The structure of $H$ is hence reminiscent of a self-similar random process
\cite{MS,ST}.
Establishing the correlation structure of the Harris paths of critical Tokunaga
processes is an interesting open problem (see Sect.\ref{open}). 
\end{Rem}

The critical Tokunaga trees introduced in Prop.~\ref{Tok_tree}
have an additional important property: the frequencies of vertex
orders in a large-order tree approximate the frequencies of orders 
in the entire space $\cT$.
To formalize this observation, let $\mu$ be the  measure on $\cT$ induced by $S^{\rm Tok}(t;c,\gamma)$, i.e. $\mu(T)={\sf P}\big(\textsc{shape}(T[S^{\rm Tok}])=T \big)$.
Next, for a fixed $K \geq 1$, let $\mu_K(T) = \mu(T|T\in\cH_K)$.
Let $V_k[K]$ denote the number of vertices of order $k \in \{1,\hdots, K\}$ in a tree generated by $\mu_K$, and let $\cV_k[K]={\sf E}(V_k[K]).$ 
Finally, we denote by $V[K]= \sum\limits_{k=1}^K V_k[K]$ the total number of non-root 
vertices, and notice that $V[K]= 2V_1[K]-1$. Thus,  $\cV[K]:={\sf E}(V[K])= 2\cV_1[K]-1$.

\begin{prop}
\label{prop:Tok1}
Let $S^{\rm Tok}(t;c,\gamma)$ be a critical Tokunaga branching process,
then
\be
\label{v1}
\lim_{ K\to\infty}\frac{\cV_k[K]}{\cV_1[K]}=2^{1-k}.
\ee
Let $T=\textsc{shape}\big(T[S^{\rm Tok}]\big)\in\cH_K$ be a tree generated by $\mu_K$, and let $v$ be a vertex selected by uniform random drawing from the
non-root vertices of $T$. Then
\be
\label{v2}
\lim_{K\to\infty}{\sf P}(v\text{ has order }k)=p_k=2^{-k}.
\ee

\end{prop}
\begin{proof}
It has been shown in \cite{KZ15t} that the mean self-similar trees satisfy
the strong Horton law:
\[\lim_{K\to\infty}\frac{\cN_k[K]}{\cN_1[K]}=R^{1-k},\text{ for any }k\ge 1.\]
Observe now that for any $T\in\cH_K$ we have
\[V_k(T) = \sum_{i=1}^{N_k(T)}(1+m_i(T)),\]
where $m_i(T)$ is the number of sub-branches that merge the $i$-th branch
of order $k \in \{1,\hdots, K\}$ in $T$, according to the proper branch labeling of 
Sect.~\ref{sec:label}.
Proposition~\ref{HBP:branch} gives
\[\cV_k[K] = \cN_k[K](1+T_1+\dots+T_{k-1}).\]
For the process $S^{\rm Tok}(t;c,\gamma)$ this implies 
\begin{eqnarray*}
\lim_{K\to\infty}\frac{\cV_k[K]}{\cV_1[K]}
&=&\lim_{K\to\infty}\frac{\cN_k[K](1+T_1\dots+T_{k-1})}{\cN_1[K]}
=R^{1-k}c^{k-1}=2^{1-k}.
\end{eqnarray*}
The statement \eqref{v2} is an immediate consequence of \eqref{v1},
since 
\[\lim_{K\to\infty}\frac{\cV_k[K]}{\cV[K]}=\lim_{K\to\infty}\frac{\cV_k[K]}{2\cV_1[K]-1}=2^{-k}\]
as $\cV_1[K] \geq 2^{K-1}$.
\end{proof}

Proposition~\ref{prop:Tok1} has an immediate extension to trees with edge lengths, which we
include here for completeness.
A tree $T\in\cL$ can be considered a metric
space with distance $d(a,b)$ between two points $a,b\in T$ defined as 
the length of the shortest path within $T$ connecting them; see \cite[Sect.~7.3]{Pitman} for details.
\begin{prop}
\label{prop:Tok2}
Consider a random tree $T=T[S^{\rm Tok}]\in\cL$ generated by a critical Tokunaga branching process $S^{\rm Tok}(t;c,\gamma)$ conditioned on the order ${\sf k}(T)=K$.
Let point $u\in T$ be sampled from a uniform density function on the metric space $T$, and let $r_u[K]$ denote the order of the edge to which the point $u$ belongs.
Then
\be
\label{p2}
\lim_{K\to\infty}{\sf P}(r_u[K]=k)=p_k=2^{-k}.
\ee
\end{prop}
\begin{proof}

Proposition~\ref{Tok_tree} establishes that the edge lengths in $T$ are i.i.d. exponential random variables. 
Thus we can generate $T$ by first sampling the combinatorial tree $\textsc{shape}(T)$ from $\cH_K$ according to
conditional measure $\mu_K(T) = \mu(T|T\in\cH_K)$, and then assigning i.i.d. exponential edge lengths. 
Provided that we already sampled $\textsc{shape}(T)$, selecting the i.i.d. edge lengths and then selecting the point $u\in T$ uniformly at random, and marking the edge that $u$ belongs to, is equivalent to selecting a random edge uniformly from the edges of $\textsc{shape}(T)$, in order of proper labeling of Sect.~\ref{sec:label}.
The order $r_u[K]$ is uniquely determined by the edge to which $u$ belongs.
The statement now follows immediately from Prop.~\ref{prop:Tok1}. 
\end{proof}

\section{Open problems}
\label{open}
We conclude with two open problems, which refer to extending
selected properties of the critical Galton-Watson tree with
independent exponential edge lengths, ${\sf GW}(0,\gamma)$,
which is a special case of the hierarchical branding process
(see Thm.~\ref{main3}), to a general case.
Our formulations are intentionally informal, reflecting multiplicity
of possible rigorous approaches to each of them. 
Here 
$S(t)=\left(\{T_k\},\{\lambda_j\},\{p_K\}\right)$
is a self-similar hierarchical branching process with
\[L=\limsup_{k\to\infty} T_k^{1/k}<\infty,\quad
p_K = p(1-p)^{K-1},\quad\lambda_j = \gamma\,\zeta^{-j}\]
for some positive $0<p<1$, $\gamma>0$, and $\zeta \ge 1\vee L$.

\begin{op}
Describe the correlation structure of the Harris path of $S(t)$.
(The critical binary Galton-Watson tree with independent exponential 
edge lengths ${\sf GW}(0,\gamma)$ corresponds to a symmetric Markov 
chain with exponential jumps $\left\{\frac{1}{2},2\,\gamma,2\,\gamma\right\}$,
see Thm.~\ref{main3}).
\end{op}

\begin{op}
Establish a proper infinite-tree limit of $S(t)$, where the edge lengths go
to zero and the tree length increases to infinity, that preserves a
suitably defined limit version of the self-similarity property.
Describe the respective limit Harris path processes.
(The Harris path of the critical binary Galton-Watson tree  ${\sf GW}(0,\gamma)$ 
can be rescaled to converge to excursion of the standard Brownian motion 
\cite{LeGall93,NP}.)
\end{op}

\section*{Acknowledgements}
We are grateful to Ed Waymire for his continuing support and encouragement. 
YK wishes to thank Tom Kurtz for providing a feedback regarding infinite dimensional population processes during  the workshop {\it Interplay of Stochastic and Deterministic Dynamics} at the Mathematical Biosciences Institute (MBI)
held February 22-26, 2016.
The authors would like to thank Jim Pitman for suggesting relevant publications. 
Comments of an anonymous referee helped us to improve presentation of the results.
We thank the participants of the conference {\it Random Trees and Maps : Probabilistic and Combinatorial Aspects}, June 6-10, 2016, Centre International de Rencontres Math\'ematiques (CIRM),
Marseille, France, and the {\it 31st Conference on Mathematical Geophysics} of the International Union
of Geodesy and Geophysics, June 6-10, 2016, Paris, France to whom we presented an early version of this work. 
The authors acknowledge financial support from the National Science Foundation, awards NSF DMS-1412557 (Y.K.) and  NSF EAR-1723033 (I.Z.)

\bibliographystyle{amsplain}

\end{document}